\documentclass[11pt]{amsart}
\usepackage{amssymb}
\usepackage{amscd}

\usepackage[all]{xy}
\input amssym.def
\input amssym
\newtheorem{theorem}{Theorem}[section]
\newtheorem{corollary}{Corollary}[section]
\newtheorem{lemma}{Lemma}[section]

\newtheorem{remark}{Remark}[section]

\newcommand{\be}{\begin{equation}}
\newcommand{\ee}{\end{equation}}

\renewcommand{\theequation}{\thesection.\arabic{equation}}
\renewcommand{\thetheorem}{\thesection.\arabic{theorem}}
\renewcommand{\theequation}{\thesection.\arabic{equation}}
\begin{document}

\title[Vertex-algebraic structure of principal subspaces]
{Vertex-algebraic structure of the principal subspaces of level one
modules for the untwisted affine Lie algebras of types $A,D,E$}

\author{C. Calinescu, J. Lepowsky and A. Milas}

\thanks{J.L. was partially supported by NSF grant DMS-0401302.  A.M.
was partially supported by NSF grant DMS-0802962.}

\begin{abstract}
Generalizing some of our earlier work, we prove natural presentations
of the principal subspaces of the level one standard modules for the
untwisted affine Lie algebras of types $A$, $D$ and $E$, and also of
certain related spaces.  As a consequence, we obtain a canonical
complete set of recursions ($q$-difference equations) for the
(multi-)graded dimensions of these spaces, and we derive their graded
dimensions.  Our methods are based on intertwining operators in vertex
operator algebra theory.
\end{abstract}

\maketitle

\renewcommand{\theequation}{\thesection.\arabic{equation}}
\renewcommand{\thetheorem}{\thesection.\arabic{theorem}}
\setcounter{equation}{0} \setcounter{theorem}{0}
\setcounter{section}{0}

\section{Introduction}

This paper is a continuation of \cite{CalLM1} and \cite{CalLM2}.  Our
aim is to continue the study of the structures called principal
subspaces, which had been introduced by Feigin-Stoyanovsky \cite{FS1},
\cite{FS2}, and related spaces, but now in the case of the basic (=
level-one standard) modules for the affine Lie algebras of types $A$,
$D$ and $E$.  Specifically, we prove a generators-and-relations result
(a presentation) for these principal subspaces, and exploit this
result to construct exact sequences and recursion relations yielding
their (multi-)graded dimensions (= generating functions of the
dimensions of the homogeneous subspaces, sometimes called
``characters'').  Our methods are based on intertwining operators in
vertex operator algebra theory.

Compared to the $\widehat{\goth{sl}(2)}$ case that was handled in
\cite{CalLM1} and \cite{CalLM2} (as well as in \cite{CLM1} and
\cite{CLM2}), the vertex-algebraic structure associated with the
principal subspaces in the present generality is more complex than the
corresponding structure in these earlier papers, and this makes the
proofs more subtle.  For a detailed introduction to and motivation of
these ideas, including historical references, we refer the reader to
\cite{CLM1}, \cite{CLM2}, \cite{Cal1}, \cite{Cal2}, \cite{CalLM1} and
\cite{CalLM2}, and Feigin-Stoyanovsky's papers \cite{FS1}, \cite{FS2}.
Presentations of principal subspaces have also been considered in
\cite{AKS}, \cite{AK} and \cite{FF}, and principal subspaces and
related structures have also been studied in \cite{P1}, \cite{G},
\cite{P2}, \cite{FFJMM} and other works.

In particular, our own interest in principal subspaces arose from the
idea developed and implemented in \cite{CLM1}, \cite{CLM2} that one
could ``explain'' the classical Rogers-Ramanujan and Rogers-Selberg
recursions for the ``sum sides'' of the Rogers-Ramanujan and
Gordon-Andrews partition identities (cf. \cite{A}) by means of exact
sequences and recursions, constructed from intertwining operators in
vertex operator algebra theory, associated with the principal
subspaces of the standard $\widehat{\goth{sl}(2)}$-modules.

Let us recall the notion of principal subspace.  Let $\goth{g}$ be a
finite-dimensional (complex) simple Lie algebra of type $A$, $D$ or
$E$.  Fix a dominant weight $\Lambda$ for $\widehat{\goth{g}}$ and
consider the standard (integrable highest weight)
$\widehat{\goth{g}}$-module $L(\Lambda)$.  Then the principal subspace
$W(\Lambda) \subset L(\Lambda)$ is defined as
\[
W(\Lambda) = U(\bar{\goth{n}}) \cdot v_{\Lambda},
\]
where $v_{\Lambda}$ is a highest weight vector of $L(\Lambda)$,
$\goth{n} \subset \goth{g}$ is the Lie subalgebra of $\goth{g}$
spanned by the root vectors for the positive roots, and
$\bar{\goth{n}}$ is the appropriate affinization of $\goth{n}$ in
$\widehat{\goth{g}}$.  Then
\[
W(\Lambda) \simeq U(\bar{\goth{n}})/{\mbox{\rm Ker}\;}f_{\Lambda},
\]
where
\[
f_{\Lambda}:U(\bar{\goth{n}}) \longrightarrow W(\Lambda)
\]
is the natural surjection that takes an element $a$ to $a \cdot
v_{\Lambda}$.

Certain well-understood formal infinite sums of elements of
$\bar{\goth{n}}$ are well known to annihilate the standard module
$L(\Lambda)$, and thus natural truncations of these formal infinite
sums lie in ${\rm Ker} \; f_{\Lambda}$.  The nontrivial part of
proving the desired presentation (the hard part) is to prove that
these truncated sums, together with obvious additional elements, {\it
generate} ${\rm Ker} \; f_{\Lambda}$.  (As is typically the case with
generators-and-relations results in mathematics, the hard part is to
prove that the ``known'' relations generate {\it all} of the relations
defining the structure being studied.)

For the case $\goth{g}=\goth{sl}(2)$, the appropriate presentations of
the principal subspaces of the standard $\widehat{\goth{g}}$-modules
of all levels were given in \cite{FS1}, \cite{FS2} (and were invoked
in the course of the proofs of the main theorems in \cite{CLM1},
\cite{CLM2}).  If one already has available an appropriate, explicit
``fermionic character formula'' for a standard module, then it is
relatively straightforward to show that the formal infinite sums
referred to above generate {\it all} the relations defining the
principal subspace, and indeed, such an explicit fermionic character
formula (as it came to be called) had been discovered in \cite{LP2},
providing enough information to justify the desired presentation of
the principal subspaces for the standard
$\widehat{\goth{sl}(2)}$-modules.  The fermionic character formula in
\cite{LP2} was in fact proved by means of the construction of what
came to be called a ``fermionic basis'' for each standard
$\widehat{\goth{sl}(2)}$-module.  These bases, which were untwisted
analogues of the fermionic bases constructed in \cite{LW1}--\cite{LW3}
(and used in those works to give vertex-operator proofs of the
Rogers-Ramanujan identities and vertex-operator interpretations of the
Gordon-Andrews-Bressoud identities), were motivated by the discovery
and use of the formal infinite sums mentioned above; these formal
infinite sums were used to contruct natural spanning sets of the
standard $\widehat{\goth{sl}(2)}$-modules (see also \cite{LP1}).  The
harder part of the work required in \cite{LP2} to complete the proof
of the main theorem constructing the fermionic bases was to prove the
{\it linear independence} of the spanning sets that had been
constructed, and as we have mentioned, the construction of these
fermionic bases (including of course the proof of their linear
independence) and the resulting fermionic character formulas allow one
to prove the nontrivial part of the expected presentation result for
the principal subspaces quite easily.  Beyond the case of
$\widehat{\goth{sl}(2)}$, it is typically difficult to construct
fermionic bases in general.

The main point our work beginning in \cite{CalLM1} and \cite{CalLM2},
and continuing in the present paper, is to prove the nontrivial part
of the presentation of the principal subspace {\it without} invoking a
theorem (or perhaps conjecture) such as a known (or perhaps proposed)
construction of a fermionic basis or a known (or proposed) fermionic
character formula.  Rather, what one really wants to do is to provide
{\it a priori} proofs of the desired presentations of the principal
subspaces, and then to use the presentation result in the course of
the construction of exact sequences and recursions whose solutions
will yield, as theorems, fermionic character formulas and fermionic
bases.  When one combines the results of \cite{CalLM1} and
\cite{CalLM2} with those of \cite{CLM1} and \cite{CLM2}, one indeed
has an a priori derivation of the desired fermionic character
formulas, without the use of explicit fermionic bases or fermionic
character formulas such as those derived in \cite{LP2}.

Specifically, the main purpose of the present paper is to give an a
priori proof of the expected presentations of the principal subspaces
of the level one standard $\widehat{\goth{g}}$-modules for types $A$,
$D$ or $E$.  Our proof is a generalization of our previous a priori
proof of this presentation for the case $\goth{g}=\goth{sl}(2)$
carried out in \cite{CalLM1}, and our methods continue the development
of the vertex-operator-algebraic ideas of \cite{CalLM1} and
\cite{CalLM2}.  Our arguments are quite delicate, and this seems to be
necessary.

When $\goth{g}$ is of type $A$, $D$ or $E$ and the standard module
$L(\Lambda)$ is of level one, the subspace $W(\Lambda)$ can be
realized inside the module $V_P$ for the lattice vertex operator
algebra $V_Q$, where $Q$ is the root lattice and $P$ is the weight
lattice of $\goth{g}$, and it is convenient for us to use this
well-known realization in this paper.

The main results of this paper are, for $\goth{g}$ of type $A$, $D$ or
$E$:

\begin{itemize}

\item[(a)] Proof of the expected presentations of the principal
subspaces for all the basic modules; that is, we explicitly describe
the left ideal ${\rm Ker} \; f_{\Lambda}$, showing that its
``obvious'' elements indeed generate it (Theorem
\ref{presentations-ps}).

\item[(b)] Construction of certain canonical exact sequences among the
principal subspaces considered in (a) (Theorem \ref{1st}).

\item[(c)] Explicit formulas for the graded dimensions of the spaces
$W(\Lambda)$ (Corollary \ref{characters}).

\item[(d)] A reformulation of (a) in terms of two-sided ideals of a
suitable completion of $U(\bar{\goth{n}})$ (Theorem
\ref{presentations-infinite}).

\end{itemize}

In fact, we in addition obtain analogues of these results for certain
spaces that we call ``principal-like'' subspaces of the $V_Q$-module
$V_P$; these subspaces arose naturally, and in fact were required, in
the course of our proofs of these results (for principal subspaces).
It is interesting that these principal-like subspaces are closely
related to versions of Kirillov-Reshetikhin modules, whose relevance
to the study of principal subspaces and so on was discussed in
\cite{AK}.

Explicit formulas such as those in (c) had already been proposed and
studied in the literature, in particular, in \cite{DKKMM},
\cite{KKMM1}, \cite{KKMM2}, \cite{KNS} and \cite{T}, in the setting of
the thermodynamic Bethe Ansatz, and formulas of this type had been
linked to principal subspaces in \cite{FS1}--\cite{FS2}, and further,
in \cite{G}, \cite{AKS}, \cite{Cal1} and \cite{FFJMM}.  Motivated by
\cite{FS1}, such formulas for the principal subspaces of certain
classes of standard modules for type $A$ were proved in \cite{G}, and
were proved by a different method more recently in \cite{Cal1} for
type $A$, level 1, a method that also included (a) and (b).  But the
main theorems in the present paper had not been proved before in the
present generality.

In the spirit of our earlier papers \cite{CalLM1}, \cite{CalLM2} (and
also in the spirit of \cite{CLM1}, \cite{CLM2}), we obtain our results
referred to in (a)--(d) in an a priori way, by combining
vertex-algebraic methods with general facts about affine Lie algebras,
and without any reference to spanning sets or bases of $W(\Lambda)$,
or to any fermionic formulas for graded dimensions.

In G. Georgiev's paper \cite{G}, which appeared shortly after
\cite{FS1}, his goal was essentially combinatorial, even though vertex
operator techniques were extensively used; combinatorial bases of
principal subspaces were obtained in the case of special families of
standard modules of type $A$.  Not surprisingly, our formulas in (c),
for type $A$, coincide with Georgiev's formulas, although our methods
are quite different.  (Georgiev was not concerned with formulating or
proving results analogous to those in (a), (b) or (d).)

The paper by E. Ardonne, R. Kedem and M. Stone \cite{AKS} gives, among
other things, a general formula for the graded dimension of
$W(\Lambda)$ for all standard modules for $\widehat{\goth{sl}(n+1)}$.
The method used to justify this formula in \cite{AKS} is reminiscent
of the method that Feigin and Stoyanovsky pursued in \cite{FS1}, where
a suitable dual of $W(\Lambda)$ was described in terms of certain
rational functions (see Theorem 3.4 of \cite{AKS}); this description
can be used to give a combinatorial interpretation of graded
dimensions.  The authors of \cite{AKS} and \cite{AK} quote a known,
natural presentation of a standard module, based on the formal
infinite sums that we have been referring to, and argue that the
subset of this set of formal infinite sums that relate in a direct way
to the principal subspace must give a complete set of defining
relations for the principal subspace (but it is a priori possible that
the other relations needed for the presentation of the standard module
could in principle contribute to further relations needed for a
correct presentation of the principal subspace, and it seems to us
that the necessity to exclude this possibility requires that a proof
of the presentation of $W(\Lambda)$ be nontrivial).  In \cite{FFJMM},
B. Feigin, E. Feigin, M. Jimbo, T. Miwa and E. Mukhin prove the
presentation of $W(\Lambda)$ for the case of $\widehat{\goth{sl}(3)}$,
all levels, in the course of establishing a formula for the graded
dimension.  In the paper \cite{FF}, B. Feigin and E. Feigin discuss
lattice vertex operator algebras, and when the Gram matrix has only
nonnegative integer entries (cf. the text before Lemma 1.1 in
\cite{FF}), a presentation result for spaces analogous to principal
subspaces is asserted, and a comparison of graded dimensions is
invoked to justify this.  Principal subspaces for general lattice
vertex operator algebras and superalgebras, and their modules and
twisted modules, are being studied in \cite{MP}.

Our paper, then, fills what appear to be some gaps in the literature
concerning principal subspaces, but what is more interesting to us is
that our methods are natural (if perhaps also rather subtle), and they
generalize considerably, as ongoing work seems to be showing.

We thank Eddy Ardonne and Rinat Kedem for helpful comments on an
earlier version of this paper.

\section{Preliminaries}

Let $\goth{g}$ be a finite-dimensional complex simple Lie algebra of
type $A$, $D$ or $E$ of rank $l$, and let $\goth{h}$ be a Cartan
subalgebra of $\goth{g}$. Let $ \{ \alpha_1, \dots, \alpha_l \}
\subset \goth{h}^{*}$ be a set of simple roots of $\goth{g}$.  Denote
by $\Delta$ the set of roots and by $\Delta_{+}$ the set of positive
roots.  We use the rescaled Killing form $\langle \cdot, \cdot
\rangle$ on $\goth{g}$ such that $\langle \alpha, \alpha \rangle=2$
for $\alpha \in \Delta$, where we identify $\goth{h}$ with
$\goth{h}^{*}$ via this form. For each root $\alpha$ fix a root vector
$x_{\alpha}$, to be rescaled later.

Denote by $\lambda_1, \dots, \lambda_l \in \goth{h} \simeq \goth{h^*}$
the corresponding fundamental weights of $\goth{g}$ (i.e., $\langle
\lambda_i, \alpha_j \rangle = \delta_{i j}$ for $i,j =1, \dots,
l$). It is also convenient to set
\[
\lambda_0=0.
\]

Let $Q=\sum_{i=1}^l \mathbb{Z} \alpha_i \subset \goth{h} \simeq
\goth{h}^{*}$ and $P=\sum_{i=1}^l \mathbb{Z} \lambda_i \subset
\goth{h} \simeq \goth{h}^{*}$ be the root and weight lattices of
$\goth{g}$, respectively.  If $\goth{g}$ of type $A_l$, $l \geq 1$, we
have $P/Q \simeq \mathbb{Z}/(l+1)\mathbb{Z}$. For $\goth{g}$ of type
$D_l$, $l \geq 4$, we have $P/Q \simeq \mathbb{Z}/4\mathbb{Z}$ or $P/Q
\simeq \mathbb{Z}/2\mathbb{Z} \times \mathbb{Z}/2\mathbb{Z}$,
depending on whether $l$ is odd or even. If $\goth{g}$ is $E_6$, $E_7$
or $E_8$ the group $P/Q$ is $\mathbb{Z}/3\mathbb{Z}$, $\mathbb{Z}/2
\mathbb{Z}$ or the trivial group, respectively.

Consider the untwisted affine Lie algebra 
\begin{equation}
\widehat{\goth{g}}= \goth{g} \otimes \mathbb{C}[t, t^{-1}] \oplus
\mathbb{C}{\bf k},
\end{equation}
where ${\bf k}$ is a nonzero central element and 
\[
[ x \otimes t^m, y \otimes t^n ] = [x, y] \otimes t^{m+n} + m\langle
x, y \rangle \delta _{m+n, 0} {\bf k}
\]
for $x, y \in \goth{g}$ and $m, n \in \mathbb{Z}$.  By adjoining the
degree operator $d$ ($[d, x \otimes t^m]=m$, $[d,{\bf k}]=0$) to the
Lie algebra $\widehat{\goth{g}}$ one obtains the affine Kac-Moody
algebra $\widetilde{\goth{g}}=\widehat{\goth{g}} \oplus \mathbb{C}d$
(cf. \cite{K}).

Set
\begin{equation}
\goth{n}=\coprod_{\alpha \in \Delta_{+}} \mathbb{C}x_{\alpha}
\end{equation}
and consider the subalgebras
\[
 \bar{\goth{n}}= \goth{n} \otimes \mathbb{C}[t, t^{-1}] ,
\]
\[
\bar{\goth{n}}_{+}=\goth{n} \otimes \mathbb{C}[t]
\]
and 
\[
\bar{\goth{n}}_{-}= \goth{n} \otimes t^{-1} \mathbb{C}[t^{-1}]
\]
of $\widehat{\goth{g}}$.  We shall frequently use the decomposition
\begin{equation}\label{decomp}
U(\bar{\goth{n}}) = U(\bar{\goth{n}}_{-}) \oplus
U(\bar{\goth{n}})\bar{\goth{n}}_{+}.
\end{equation}
The affine Lie algebra $\widehat{\goth{g}}$ has the decomposition
\[
\widehat{\goth{g}}= \widehat{\goth{g}}_{<0} \oplus
\widehat{\goth{g}}_{\geq 0},
\]
where
\[
\widehat{\goth{g}}_{<0}= \goth{g} \otimes t^{-1} \mathbb{C} [t^{-1}]
\]
and
\[ \widehat{\goth{g}}_{\geq 0}= \goth{g} \otimes \mathbb{C}[t] \oplus
\mathbb{C}{\textbf k}.
\]

The form $\langle \cdot, \cdot \rangle$ on $\goth{g}$ extends
naturally to $\goth{h} \oplus \mathbb{C}{\bf k} \oplus \mathbb{C}d$,
with $\langle {\bf k}, d \rangle = 1$.  We shall identify $\goth{h}
\oplus \mathbb{C}{\bf k} \oplus \mathbb{C}d$ with its dual space
$(\goth{h} \oplus \mathbb{C}{\bf k} \oplus \mathbb{C}d)^{*}$ via this
form. As usual, we denote by $ \alpha_0, \alpha_1, \dots, \alpha_l \in
(\goth{h} \oplus \mathbb{C}{\bf k} \oplus \mathbb{C}d)^{*}$ the
corresponding simple roots, and by $\Lambda_0, \Lambda_1, \dots,
\Lambda_l \in (\goth{h} \oplus \mathbb{C}{\bf k} \oplus
\mathbb{C}d)^{*}$ the corresponding fundamental weights.  Then
$\langle \Lambda_0, {\bf k} \rangle =1$; for $i=1, \dots, l$, $\langle
\Lambda_i, {\bf k} \rangle =k_i$, where $k_i \geq 1$ is the
coefficient of $\alpha_i$ in the expansion of the highest root;
$\langle \Lambda_i, \alpha_j \rangle= \delta_{ij}$ for $i, j =0,
\dots, l$; and
\begin{equation} \label{d}
\langle \Lambda_0, d \rangle =0, \; \; \; 
\langle \Lambda_i, d \rangle
=-\frac{1}{2} \langle \Lambda_i, \Lambda_i \rangle \; \;  \mbox{for}
\; \; i=1, \dots,
l.
\end{equation}
Moreover, for $i=0, \dots, l$,
\[
\Lambda_i=\lambda_i+\langle \Lambda_i, d \rangle {\bf k}+ \langle
\Lambda_i, {\bf k} \rangle d,
\]
so that in particular, $\Lambda_0=d$.  In this paper we shall focus on
the level one standard $\widehat{\goth{g}}$-modules (also called the
basic $\widehat{\goth{g}}$-modules), denoted $L(\Lambda)$, where
$\Lambda$ is one of the fundamental weights $\Lambda_i$ such that
$\langle \Lambda, {\bf k} \rangle=1$. Let $v_{\Lambda}$ be a highest
weight vector of such a module $L(\Lambda)$.  We shall normalize
$v_{\Lambda}$ later. We view each of the basic
$\widehat{\goth{g}}$-modules as a $\widetilde{\goth{g}}$-module, where
$d$ acts according to (\ref{d}) on $v_{\Lambda}$. The order of $P/Q$
gives the number of inequivalent standard $\widehat{\goth{g}}$-modules
of level one: For $\goth{g}=A_l$ there are $l+1$ level one standard
$\widehat{\goth{g}}$-modules: $L(\Lambda_0), L(\Lambda_1), \dots ,
L(\Lambda_l)$. For $\goth{g}=D_l$ there are four such modules:
$L(\Lambda_0)$, $L(\Lambda_1)$, $L(\Lambda_{l-1})$ and
$L(\Lambda_l)$. If $\goth{g}$ is $E_6$ there are three:
$L(\Lambda_0)$, $L(\Lambda_1)$ and $L(\Lambda_6)$. There are two basic
$\widehat{\goth{g}}$-modules, $L(\Lambda_0)$ and $L(\Lambda_1)$, when
$\goth{g}$ is $E_7$, and only one such module, $L(\Lambda_0)$, if
$\goth{g}$ is $E_8$.  (Cf. \cite{K}.)

For each fundamental weight $\Lambda_i$ of $\widehat{\goth{g}}$ such
that $\langle \Lambda_i, {\bf k} \rangle=1$ consider the principal
subspace $W(\Lambda_i)$ of the level one standard module
$L(\Lambda_i)$ in the sense of \cite{FS1}--\cite{FS2} (see also
\cite{CalLM2} for a natural generalization of this notion):
\begin{equation}
W(\Lambda_i)= U(\bar{\goth{n}}) \cdot v_{\Lambda_i}.
\end{equation}
By the highest weight vector property we have
\begin{equation}
W(\Lambda_i)= U(\bar{\goth{n}}_{-}) \cdot v_{\Lambda_i}.
\end{equation}
As in \cite{CalLM1}--\cite{CalLM2} we consider the surjective maps
\begin{eqnarray} \label{surj1}
F_{\Lambda_{i}}: U(\widehat{\goth{g}}) & \longrightarrow &
L(\Lambda_{i}) \\ a &\mapsto& a \cdot v_{\Lambda_{i}}. \nonumber
\end{eqnarray}
Restrict $F_{\Lambda_{i}}$ to $U(\bar{\goth{n}})$ and denote these
(surjective) restrictions by $f_{\Lambda_{i}}$ :
\begin{eqnarray} \label{surj2}
f_{\Lambda_{i}}: U(\bar{\goth{n}}) & \longrightarrow &
W(\Lambda_{i})\\ a & \mapsto & a \cdot v_{\Lambda_{i}}. \nonumber
\end{eqnarray}
In this paper we will give a precise description of the
kernels $\mbox{Ker} \; f_{\Lambda_{i}}$, and thus a presentation of
the principal subspaces $W(\Lambda_i)$.

\begin{remark}
\rm In \cite{CalLM1} and \cite{CalLM2}, we used the symbol
$f_{\Lambda}$ for the restriction of $F_{\Lambda}$ to
$U(\bar{\goth{n}}_{-})$ rather than to $U(\bar{\goth{n}})$, but since
$\bar{\goth{n}}$ is no longer abelian in general, the present maps
$f_{\Lambda}$ are the appropriate ones.
\end{remark}

As in \cite{CalLM1}--\cite{CalLM2} we will sometimes use generalized
Verma modules for $\widehat{\goth{g}}$, in the sense of \cite{L1},
\cite{GL} and \cite{L2}, in our formulations and proofs.  The
generalized Verma module $N(\Lambda_0)$ is defined as the induced
$\widehat{\goth{g}}$-module
\begin{equation}
N(\Lambda_0)=U(\widehat{\goth{g}}) \otimes_{U(\widehat{\goth{g}}_{\geq
0})}\mathbb{C}v_{\Lambda_0}^N,
\end{equation}
where $\goth{g} \otimes \mathbb{C}[t]$ acts trivially and ${\bf k}$
acts as the scalar $1$ on $\mathbb{C}v_{\Lambda_0}^N$;
$v_{\Lambda_0}^N$ is a highest weight vector.  More generally, for any
fundamental weight $\Lambda_i$ such that $\langle \Lambda_i, {\bf k}
\rangle =1$, define the generalized Verma module
\begin{equation}
N(\Lambda_i)= U(\widehat{\goth{g}})
\otimes_{U(\widehat{\goth{g}}_{\geq 0})} U_i
\simeq U(\widehat{\goth{g}}_{< 0}) \otimes U_i,
\end{equation}
where $U_i$ is a copy of the finite-dimensional irreducible
$\goth{g}$-module $U(\goth{g}) \cdot v_{\Lambda_{i}} \subset
L(\Lambda_{i}),$ with highest weight vector now called
$v_{\Lambda_i}^N$, and where $\goth{g} \otimes t \mathbb{C} [t]$ acts
trivially and ${\bf k}$ by $1$.  We view all these generalized Verma
modules as $\widetilde{\goth{g}}$-modules, where $d$ acts according to
(\ref{d}) on $v_{\Lambda_i}^N$.  Continuing to generalize
\cite{CalLM1}--\cite{CalLM2}, we introduce
\begin{equation}\label{WN}
W^N(\Lambda_i)= U(\bar{\goth{n}}) \cdot v_{\Lambda_i}^N
\simeq U(\bar{\goth{n}}_{-}) \cdot v_{\Lambda_i}^N,
\end{equation}
the principal subspace of $N(\Lambda_i)$; this is naturally an
$\bar{\goth{n}} \oplus {\goth{h}} \oplus \mathbb{C}{\bf k} \oplus
\mathbb{C}d$-submodule of $N(\Lambda_i)$.  We have the natural
surjective $\widehat{\goth{g}}$-module maps
\begin{eqnarray} 
F^N_{\Lambda_{i}}: U(\widehat{\goth{g}}) & \longrightarrow &
N(\Lambda_{i}) \\ a &\mapsto& a \cdot v^N_{\Lambda_{i}} \nonumber
\end{eqnarray}
and their restrictions to $U(\bar{\goth{n}})$,
\begin{eqnarray} \label{f_N_map}
f^N_{\Lambda_{i}}: U(\bar{\goth{n}}) & \longrightarrow &
W^N(\Lambda_{i})\\ a & \mapsto & a \cdot v^N_{\Lambda_{i}}, \nonumber
\end{eqnarray}
where $\langle \Lambda_i, {\bf k} \rangle =1$.  For such fundamental
weights $\Lambda_i$ we have the natural surjective
$\widetilde{\goth{g}}$-module maps
\begin{eqnarray} 
\Pi_{\Lambda_i}: N(\Lambda_i)& \longrightarrow & L(\Lambda_i) \\ a
\cdot v^N_{\Lambda_i} &\mapsto& a \cdot v_{\Lambda_{i}} \nonumber
\end{eqnarray}
for $a \in U(\widehat{\goth{g}})$ and their restrictions to
$W^N(\Lambda_i)$,
\begin{eqnarray}  \label{pi_map}
\pi_{\Lambda_i}: W^N(\Lambda_i)& \longrightarrow &
W(\Lambda_i);
\end{eqnarray}
these are $U(\bar{\goth{n}} \oplus {\goth{h}} \oplus \mathbb{C}{\bf k}
\oplus \mathbb{C}d)$-module surjections.

Throughout this paper we will write $x(m)$ for the action of $x
\otimes t^m \in \widehat{\goth{g}}$ on any
$\widehat{\goth{g}}$-module, for $x \in \goth{g}$ and $m \in
\mathbb{Z}$. In particular, we have the operator $x_{\alpha}(m)$, the
image of $x_{\alpha} \otimes t^m$, for $\alpha \in \Delta$.  Sometimes
we will write $x(m)$ for the Lie algebra element $x \otimes t^m$
itself; it will be clear {}from the context whether $x(m)$ is an
operator or a Lie algebra element.

Now we recall the constructions of lattice vertex operator algebras
and their modules from Section 7.1 and Chapter 8 of \cite{FLM2}
(cf. Sections 6.4 and 6.5 of \cite{LL}).  We work in the setting of
\cite{LL}, adapted to our situation.  Consider
$\widehat{\goth{h}}=\goth{h} \otimes \mathbb{C}[t, t^{-1}] \oplus
\mathbb{C}{\bf k}$ and its irreducible induced module
\begin{equation}
M(1)=U(\widehat{\goth{h}}) \otimes_{U(\goth{h} \otimes \mathbb{C}[t]
\oplus \mathbb{C} {\bf k})} \mathbb{C},
\end{equation}
where $\goth{h} \otimes \mathbb{C}[t]$ acts trivially and ${\bf k}$
acts as 1 on the one-dimensional module $\mathbb{C}$. The space $M(1)$
can be identified with the symmetric algebra
$S(\widehat{\goth{h}}_{-})$, where
\[
\widehat{\goth{h}}_{-}=\goth{h} \otimes t^{-1} \mathbb{C}[t^{-1}].
\]

We shall fix $s > 0$ and a central extension $\widehat{P}$ of the
weight lattice $P$ (and by restriction this gives a central extension
of the root lattice $Q$) by the finite cyclic group $\langle \kappa
\rangle= \langle \kappa \; | \: {\kappa}^s=1 \rangle$ of order $s$,
\[
1 \rightarrow \langle \kappa \rangle \rightarrow \widehat{P}
\bar{\longrightarrow} P \rightarrow 1,
\]
satisfying the condition (\ref{com}) below.  Let $c_0: P \times P
\longrightarrow \mathbb{Z}/s\mathbb{Z}$ be the associated commutator
map, so that $aba^{-1}b^{-1}={\kappa}^{c_{0}(\bar{a},\bar{b})}$ for
$a,b \in \widehat{P}$, let $\nu_s \in \mathbb{C}^{\times}$ be a
primitive $s^{\rm th}$ root of unity, and define the map $c:P \times P
\longrightarrow \mathbb{C}^{\times}$ by $c(\alpha,
\beta)=\nu_s^{c_0(\alpha, \beta)}$ for $\alpha, \beta \in P$.  Then
$c(\cdot,\cdot)$ is an alternating $\mathbb{Z}$-bilinear map from $P
\times P$ to the multiplicative group $\mathbb{C}^{\times}$.  We
assume the condition
\begin{equation} \label{com}
c(\alpha, \beta)=(-1)^{\langle \alpha, \beta \rangle} \; \; \;
\mbox{for} \; \; \; \alpha, \beta \in Q;
\end{equation}
there indeed exists $s > 0$ together with such a central extension
$\widehat{P}$ (see Remark 6.4.12 in \cite{LL}).

Define the faithful character $\chi: \langle \kappa \rangle
\longrightarrow \mathbb{C}^{\times}$ by the condition $\chi(\kappa)=
\nu_s$, and denote by $\mathbb{C} \{ P \}$ the induced
$\widehat{P}$-module
$\mathbb{C}[\widehat{P}]\otimes_{\mathbb{C}[\langle \kappa
\rangle]}\mathbb{C}_{\chi}$, where $\mathbb{C}_{\chi}$ is the
one-dimensional space $\mathbb{C}$ viewed as a $\langle \kappa
\rangle$-module (i.e., $ \kappa \cdot 1=\nu_s$).  Then the space
\[
V_Q=M(1) \otimes \mathbb{C}\{ Q \}
\]
carries a natural vertex operator algebra structure, with
$1$ as vacuum vector, and the space 
\[
V_P=M(1) \otimes \mathbb{C}\{ P \}
\]
is a $V_Q$-module in a natural way, as specified in \cite{FLM2} and
Sections 6.4 and 6.5 of \cite{LL}, in particular, Theorems 6.5.1,
6.5.3 and 6.5.20 of \cite{LL}.  For the ``well known,'' but
nontrivial, natural uniqueness of these structures of vertex operator
algebra and module, see Remark 6.5.4, Proposition 6.5.5, Remark 6.5.6
and Remark 6.5.25 of \cite{LL}.

We recall certain features of this structure from \cite{LL}.  For
convenience, choose a section
\begin{eqnarray} \label{section}
e: P & \longrightarrow & \widehat{P} \\ \nonumber
\alpha & \mapsto & e_{\alpha}, \nonumber
\end{eqnarray}
normalized by the condition $e_0 = 1$, and denote by $\epsilon_0: P
\times P \longrightarrow \mathbb{Z}/s\mathbb{Z}$ the corresponding
$2$-cocycle, defined by the condition
$e_{\alpha}e_{\beta}={\kappa}^{\epsilon_0(\alpha, \beta)}
e_{\alpha+\beta}$ for $\alpha, \beta \in P $.  Define $\epsilon:P
\times P \longrightarrow \mathbb{C}^{\times}$ by $\epsilon(\alpha,
\beta)=\nu_s^{\epsilon_0(\alpha, \beta)}$.  Then for any $\alpha,
\beta \in P$ we have
\begin{equation} \label{e-c}
\epsilon(\alpha, \beta)  / \epsilon(\beta, \alpha)= c(\alpha, \beta)
\end{equation}
and 
\begin{equation} \label{epsilon-zero}
\epsilon(\alpha, 0)=\epsilon(0, \alpha)=1.
\end{equation}
The choice of the
section (\ref{section}) allows us to identify $\mathbb{C}\{P \}$ with
the group algebra $\mathbb{C}[P]$, viewed as a vector space, by the
linear isomorphism
\begin{eqnarray} \label{identification}
\mathbb{C}[P] & \longrightarrow & \mathbb{C}\{ P \} \\ 
e^{\alpha} & \mapsto & \iota(e_{\alpha}) \nonumber
\end{eqnarray}
for $\alpha \in P$, where, for $a \in \widehat{P}$, we set $\iota(a)=
a \otimes 1 \in \mathbb{C} \{ P \}$.  The action of $\widehat{P}$ on
$\mathbb{C}[P]$ is given by $e_{\alpha} \cdot e^{\beta}=
{\epsilon(\alpha, \beta)}e^{\alpha + \beta}$, $\kappa \cdot e^{\beta}=
\nu_{s} e^{\beta}$ for $\alpha, \beta \in P$, and as operators on
$\mathbb{C}[P] \simeq \mathbb{C}\{P \}$ we have
\begin{equation} \label{e-multiplication}
e_{\alpha}e_{\beta}=\epsilon(\alpha, \beta)e_{\alpha+\beta}.
\end{equation}  
We
also have the identification $\mathbb{C}[Q] \simeq \mathbb{C}\{ Q \}$
and the identifications
\[
V_Q=M(1) \otimes \mathbb{C}[Q]
\]
and
\[
V_P=M(1) \otimes \mathbb{C}[P].
\]
For $h \in \goth{h}$ and $n \in \mathbb{Z}$, we have the standard
operators $h(n)$ on $V_P$ (recall formulas (6.4.47) and (6.4.48) in
\cite{LL}), providing $V_P$ with $\widehat{\goth{h}}$-module
structure. The module $V_P$ for the vertex operator algebra $V_Q$ is
the direct sum of the irreducible $V_Q$-modules $M(1) \otimes
\mathbb{C}[Q]e^{\lambda_i}$, where $i$ ranges through the indices such
that $\langle \Lambda_i, {\bf k} \rangle =1$. For all $i=0, \dots, l$,
including those for which $\langle \Lambda_i, {\bf k} \rangle=1$, set
\begin{equation}\label{VQ}
V_Qe^{\lambda_i}=M(1) \otimes \mathbb{C}[Q]e^{\lambda_i}.
\end{equation}

For $\lambda \in P$ we have the vertex operator 
\begin{equation} \label{vertex}
Y(\iota(e_{\lambda}), x)=E^{-}(-\lambda, x)E^{+}(-\lambda,
x)e_{\lambda}x^{\lambda}
\end{equation}
(see formula (6.4.65) in \cite{LL}), where
\[
E^{\pm}(-\lambda, x)=\mbox{exp} \left ( \sum_{\pm n >0}
\frac{-\lambda(n)}{n} x^{-n} \right ) \in( \mbox{End} \; V_P) [[x,
x^{-1}]]
\]
and the operator $x^{\lambda}$ is defined by 
$$
x^{\lambda} (v \otimes \iota(e_{\beta}))=x^{\langle \lambda, \beta
\rangle} (v \otimes \iota(e_{\beta}))
$$ 
for $v \in M(1)$ and $\beta \in P$.  Using the identification
(\ref{identification}) we we shall write $Y(e^{\lambda}, x)$ instead
of $Y(\iota(e_{\lambda}), x)$, for convenience. In particular, for any
$\alpha \in \Delta$ we have the operators $x_{\alpha}(m)$ defined by
\begin{equation} 
Y(e^{\alpha}, x)=
\sum_{m \in \mathbb{Z}} x_{\alpha}(m) x^{-m-1}.
\end{equation}
These operators together with the action of $\widehat{\goth{h}}$ give
$V_P$ a $\widehat{\goth{g}}$-module structure, and we identify
$x_{\alpha}(0)$ with the root vector $x_{\alpha} \in \goth{g}$.
Recall from Proposition 6.4.5 of \cite{LL} that
\begin{equation}
x^{\lambda}e_{\beta}=x^{\langle \lambda, \beta \rangle}
e_{\beta}x^{\lambda}
\end{equation}
and 
\begin{equation} \label{lambda-e}
\lambda(m) e_{\beta}=e_{\beta} \lambda(m)
\end{equation}
for all $\lambda, \beta \in P$ and $m \in \mathbb{Z}$. Using
(\ref{e-c}), (\ref{e-multiplication}) and
(\ref{vertex})-(\ref{lambda-e}) we obtain, for $\alpha \in \Delta$,
\begin{equation}\label{almost-comm}
x_{\alpha}(m)e_{\beta}=c(\alpha, \beta) e_{\beta}x_{\alpha}(m+\langle
\alpha, \beta \rangle).
\end{equation} 

We take
\[\omega=\frac{1}{2 } \sum_{i=1}^{l}  u^{(i)}(-1)^2 {\bf 1}\]
for the standard conformal vector, where $\{u^{(i)}, \dots, u^{(l)}
\}$ is an orthonormal basis of $\goth{h}$ (recall formula (6.4.9) in
\cite{LL}), so that the operators $L(n)$ defined by
\begin{equation}\label{Yomega}
Y(\omega,x)=\sum_{n \in \mathbb{Z}} L(n)x^{-n-2}
\end{equation}
provide a representation of the Virasoro algebra of central charge
$l$.

For each fundamental weight $\Lambda_i$ with $\langle \Lambda_i, {\bf
k} \rangle =1$ we may and do identify the level one standard
$\widehat{\goth{g}}$-module $L(\Lambda_i)$ with the irreducible
$V_Q$-submodule $V_Qe^{\lambda_i}$ of $V_P$ (recall (\ref{VQ})), so
that in particular $L(\Lambda_0)=V_Q$, and we take as its highest
weight vector
\begin{equation}
v_{\Lambda_i}=e^{\lambda_i};
\end{equation}
in particular,
\[
v_{\Lambda_0}=1.
\]
For $\goth{g}=A_l$, the irreducible $L(\Lambda_0)$-modules (up to
isomorphism) are $L(\Lambda_0), \dots, L(\Lambda_l)$. For
$\goth{g}=D_l$ the spaces $L(\Lambda_0)$, $L(\Lambda_1)$,
$L(\Lambda_{l-1})$ and $L(\Lambda_l)$ are the irreducible
$L(\Lambda_0)$-modules. For $\goth{g}=E_6$, $L(\Lambda_0)$,
$L(\Lambda_1)$ and $L(\Lambda_6)$ are the irreducible
$L(\Lambda_0)$-modules. When $\goth{g}=E_7$, $L(\Lambda_0)$ and
$L(\Lambda_1)$ are the irreducible $L(\Lambda_0)$-modules, and for
$\goth{g}=E_8$, $L(\Lambda_0)$ is the only irreducible
$L(\Lambda_0)$-module. (See \cite{D}, \cite{DL}, \cite{DLM} and
\cite{LL}.)

The generalized Verma module $N(\Lambda_0)$ carries a natural
structure of vertex operator algebra with $v_{\Lambda_0}^N$ as vacuum
vector and with central charge $l$ (see Theorem 6.2.18 in \cite{LL}),
and for each fundamental weight $\Lambda_i$ such that $\langle
\Lambda_i, {\bf k} \rangle=1$ the spaces $N(\Lambda_i)$ are naturally
modules for $N(\Lambda_0)$ (see Theorem 6.2.21 of \cite{LL}).

If $\langle \Lambda_i, {\bf k} \rangle >1$, then $e^{\lambda_i} \in
V_P$ is not a highest weight vector for $\widehat{\goth{g}}$, but it
is a highest weight vector for $\bar{\goth{n}}$ in the sense that it
is annihilated by $\bar{\goth{n}}_{+}$ and its span is preserved by
${\goth{h}} \oplus \mathbb{C}{\bf k} \oplus \mathbb{C}d$.  We shall
use the notation
\begin{equation} \label{v_lambda_i}
v_{\lambda_i}=e^{\lambda_i} 
\end{equation}
for $i=0, \dots, l$, so that $v_{\lambda_i}$ agrees with the highest
weight vector $v_{\Lambda_i}$ if $\langle \Lambda_i, {\bf k} \rangle
=1$.  We now generalize the notion of principal subspace as follows:
We define
\begin{equation}
W(\lambda_i)=U(\bar{\goth{n}}) \cdot v_{\lambda_i} \subset V_P
\end{equation}
for each $i=0, \dots, l$, and we call these the {\it principal-like
subspaces}. More generally, for any $\lambda \in P$ we have the
{\it principal-like subspace}
\begin{equation}
W(\lambda)=U(\bar{\goth{n}}) \cdot e^{\lambda} \subset V_P. 
\end{equation}
We also have the {\it principal-like subspace}
\begin{equation}
W^N(\lambda_i)=U(\bar{\goth{n}}) \otimes_{U(\bar{\goth{n}}_{+})}
\mathbb{C} v_{\lambda_i}^N \simeq U(\bar{\goth{n}}_{-}) \cdot
v_{\lambda_i}^N,
\end{equation}
where $v_{\lambda_i}^N$ is a highest weight vector for
$\bar{\goth{n}}$ and $i=0, \dots, l$; if $\langle \Lambda_i, {\bf k}
\rangle =1$ we take $v_{\lambda_i}^N$ to be the vector
$v_{\Lambda_i}^N$ used in (\ref{WN}).  We view $W^N(\Lambda_i)$ as an
$\bar{\goth{n}} \oplus {\goth{h}} \oplus \mathbb{C}{\bf k} \oplus
\mathbb{C}d$-module, as in (\ref{WN}), where $\goth{h}$, ${\bf k}$ and
$d$ act on $v^N_{\Lambda_i}$ as they do on $e^{\lambda_i}$.

\begin{remark}  \label{w=w}
\rm For each fundamental weight $\Lambda_i$ with $\langle \Lambda_i,
{\bf k} \rangle =1$ the principal-like subspace $W(\lambda_i)$ agrees
with the principal subspace $W(\Lambda_i)$, and $W^N(\lambda_i)$
agrees with $W^N(\Lambda_i)$.
\end{remark}

Standard arguments show that $W(\lambda_0)$ is a vertex subalgebra of
$L(\Lambda_0)$ and that each $W(\lambda_i)$, $i=0, \dots, l$, is a
module for this vertex algebra.  Moreover, each $W(\lambda_i)$ is
preserved by $L(0)$ and by the action of $\goth{h}$.

Generalizing (\ref{surj2}), (\ref{f_N_map}) and (\ref{pi_map}) we have
the natural maps
\begin{eqnarray} 
f_{\lambda_{i}}: U(\bar{\goth{n}}) & \longrightarrow &
W(\lambda_{i})\\ a & \mapsto & a \cdot v_{\lambda_i}, \nonumber
\end{eqnarray}
\begin{eqnarray} 
f^N_{\lambda_{i}}: U(\bar{\goth{n}}) & \longrightarrow &
W^N(\lambda_{i})\\ a & \mapsto & a \cdot v_{\lambda_i}^N, \nonumber
\end{eqnarray}
\begin{eqnarray} 
\pi_{\lambda_{i}}: W^N(\lambda_i) & \longrightarrow & W(\lambda_{i})\\
a \cdot v_{\lambda_i}^N & \mapsto & a \cdot v_{\lambda_i}, \nonumber
\end{eqnarray}
where $a \in U(\bar{\goth{n}})$, for $i=0, \dots, l$.  
\begin{remark} \label{f=f}
\rm Note that the maps $f_{\lambda_i}$, $f^N_{\lambda_i}$ and
$\pi_{\lambda_i}$ indeed agree with the maps $f_{\Lambda_i}$,
$f^N_{\Lambda_i}$ and $\pi_{\Lambda_i}$, respectively, when $\langle
\Lambda_i, {\bf k} \rangle =1$.
\end{remark}

Recall from formula (5.1.5), Remark 5.4.2 and Proposition 5.4.7 of
\cite{FHL} the (nonzero) intertwining operator
\begin{equation} \label{intertwining_op}
{\mathcal Y}(\cdot , x): V_P \longrightarrow \mbox{Hom} \;
(V_Q, V_P)[[x,x^{-1}]]
\end{equation}
defined by
\begin{equation}
{\mathcal Y} (w, x)v=e^{xL(-1)}Y(v, -x)w, \; \; \; w \in V_P, \; \; v
\in V_Q,
\end{equation}
where $L(-1)$ is the usual Virasoro algebra operator (recall
(\ref{Yomega})).  In particular, we have
\begin{equation} \label{inter}
{\mathcal Y}(e^{\lambda_i}, x): V_Q \longrightarrow V_P ((x))
\end{equation}
for $i=0, \dots, l$.  By a standard argument (cf. \cite{CLM1}),
\[
[x_{\alpha}(m),\mathcal{Y}(e^{\lambda_i},x)]=0
\]
for all $i$. We denote the constant term (the coefficient of
$x^0$) of ${\mathcal Y}(e^{\lambda_i}, x)$ by ${\mathcal
Y}_c(e^{\lambda_i}, x)$. Then we have a surjection
\begin{equation} \label{intertwining}
{\mathcal Y}_c(e^{\lambda_i}, x): W(\lambda_0) \longrightarrow
W(\lambda_i),
\end{equation}
since this map sends $v_{\lambda_0}$ to $v_{\lambda_i}$ and ${\mathcal
Y}_c(e^{\lambda_i}, x)$ commutes with the action of $\bar{\goth{n}}$.
Thus
\begin{equation} \label{inclusion_kernels}
\mbox{Ker} \; f_{\lambda_0} \subset \mbox{Ker} \; f_{\lambda_i}.
\end{equation}
Indeed, let $a \in U(\bar{\goth{n}})$ such that $a \in \mbox{Ker} \;
f_{\lambda_0}$; then $a \cdot e^{\lambda_0}=0$. By applying the map
(\ref{intertwining}) and using its properties we have $a \cdot
e^{\lambda_i}=0$, and thus $a \in \mbox{Ker} \; f_{\lambda_i}$.

\begin{remark} \rm The $V_Q$-module $V_P$ has a structure of an
abelian intertwining algebra, as defined and described in Chapter 12
of \cite{DL}. The intertwining operators (\ref{intertwining_op}) and
(\ref{inter}) are operators in this abelian intertwining algebra, or
more precisely, the restrictions of these intertwining operators to
individual sectors corresponding to the cosets of $Q$ in $P$ are the
operators that are part of the abelian intertwining algebra, modulo
certain normalizations.
\end{remark}

The vector space $V_P$ has a natural grading defined by the action of
the standard Virasoro algebra operator $L(0)$ introduced earlier,
referred to as the grading by {\it weight}. In particular, we have
\[
\mbox{wt} \; e^{\lambda}= \frac{1}{2} \langle \lambda, \lambda \rangle
\]
for any $\lambda \in P$, and
\[
\mbox{wt} \; x_{\alpha}(m)=-m
\]
for any $\alpha \in \Delta$ and $m \in \mathbb{Z}$.  The
$L(0)$-eigenspaces of $V_P$ coincide with the eigenspaces for the
negative $-d$ of the degree operator $d$ (with the same respective
eigenvalues).  There are also $l$ gradings by {\it charge} on $V_P$,
given by the eigenvalues of the operators $\lambda_i=\lambda_i(0)$,
$i=1, \dots, l$. We will refer to these as the $\lambda_i$-{\it
charge} gradings. The weight and charge gradings are compatible.  For
any $m \in \mathbb{Z}$, $x_{\alpha_i}(m)$, viewed as either an
operator or as an element of $U(\bar{\goth{n}})$, has weight $-m$ and
charge $\delta_{ij}$ with respect to $\lambda_j$.  Also,
$e^{\lambda_i}$ has charge $\langle \lambda_j, \lambda_i \rangle$ with
respect to $\lambda_j$.  We define {\it total charge} as the sum of
the $\lambda_i$-charges.  We restrict these gradings to the
principal-like subspaces $W(\lambda_i)$ for $i=0, \dots, l$, and in
particular to the principal subspaces $W(\Lambda_i)$, and we give the
spaces $W^N(\lambda_i)$ and $W^N(\Lambda_i)$ the analogous gradings.

\begin{remark} \label{L(0)}
\em Just as in \cite{CalLM1}--\cite{CalLM2}, we have that $\mbox{Ker}
\; f_{\Lambda_i}$ and $\mbox{Ker} \; \pi_{\Lambda_i}$ are graded by
weight and by $\lambda_j$-charge for $j=1, \dots, l$, and these
gradings are compatible. More generally, this assertion holds for
$\mbox{Ker} \; f_{\lambda_i}$ and $\mbox{Ker} \; \pi_{\lambda_i}$ for
each $i=0, \dots, l$.
\end{remark}

\section{Ideals and morphisms}
\setcounter{equation}{0}
 
We consider the following formal infinite sums of operators:
\begin{equation}  \label{R}
R_{t}^i= \sum_{m_1+m_2=-t}x_{\alpha_i}(m_1) x_{\alpha_i}(m_2), \ \ t
\in \mathbb{Z}, \ \ i=1,\dots,l;
\end{equation}
each $R_t^i$ acts naturally on any highest weight
$\widehat{\goth{g}}$-module and more generally on any
$\bar{\goth{n}}$-module on which the formal sum terminates when
applied to any vector. We truncate $R_t^i$ as follows:
\begin{equation} \label{r}
R_{-1, t}^i= \sum_{m_1, m_2 \leq -1, \; m_1+m_2=-t}x_{\alpha_i}(m_1)
x_{\alpha_i}(m_2)
\end{equation}
for $t \in \mathbb{Z}$ (so that $R_{-1, t}^i=0$ unless $t \geq 2$) and
$i=1, \dots , l$. When $\goth{g}=\goth{sl}(2)$, and thus only $i=1$ is
relevant, these are the formal sums introduced in \cite{CalLM1} and
denoted by $R_t^0$.  We shall often be viewing $R_{-1,t}^i$ as an
element of $U(\bar{\goth{n}})$, and in fact of
$U(\bar{\goth{n}}_{-})$, rather than as an endomorphism of a
$\widehat{\goth{g}}$-module.  Note that $R_{t}^i$ and $R_{-1, t}^i$
have weight $t$, charge 2 with respect to $\lambda_i$, and charge $0$
with respect to $\lambda_j$ for $j \neq i$.  In this paper we will
give two statements for presentations of the principal subspaces
$W(\Lambda_i)$, where $\Lambda_i$ is a fundamental weight such that
$\langle \Lambda_i, {\bf k} \rangle=1$, and, more generally, of the
principal-like subspaces $W(\lambda_i)$, $0 \leq i \leq l$. One
statement involves left ideals of $U(\bar{\goth{n}})$ generated by
elements of type (\ref{r}), while the other statement uses two-sided
ideals of a certain completion of $U(\bar{\goth{n}})$, ideals
generated by the formal infinite sums (\ref{R}).

Let $J$ be the left ideal of $U(\bar{\goth{n}})$ generated by the
elements $R_{-1, t}^i$ for $t \geq 2$ and $i=1,\dots,l$:
\begin{equation}
J = \sum_{i=1}^l \sum_{t \geq 2} U(\bar{\goth{n}}) R_{-1,t}^i.
\end{equation}
By analogy with the corresponding constructions in
\cite{CalLM1}--\cite{CalLM2} (which involved $U(\bar{\goth{n}}_-)$),
we set
\begin{equation} \label{I_0}
I_{\lambda_{0}}= J+U(\bar{\goth{n}})\bar{\goth{n}}_{+} \subset
U(\bar{\goth{n}})
\end{equation}
and
\begin{equation} \label{I_lambda_i}
I_{\lambda_i}=J+U(\bar{\goth{n}})\bar{\goth{n}}_{+}
+U(\bar{\goth{n}})x_{\alpha_i}(-1)
=I_{\lambda_{0}} + U(\bar{\goth{n}})x_{\alpha_i}(-1)
\subset U(\bar{\goth{n}}), \; \; i=1, \dots, l.
\end{equation}

\begin{remark}
\rm Since in this paper we are concerned with more general structures
than the principal subspaces of the level one standard
$\widehat{\goth{g}}$-modules, and since these structures,
$W(\lambda_i)$, are indexed by the fundamental weights of $\goth{g}$
rather than the fundamental weights of $\widehat{\goth{g}}$, we use
the notation $I_{\lambda_i}$ rather than $I_{\Lambda_i}$. Note that
the ideals $I_{\lambda_0}$ and $I_{\lambda_1}$ are the same as
$I_{\Lambda_0}$ and $I_{\Lambda_1}$ in \cite{CalLM1} modulo
$U(\bar{\goth{n}})\bar{\goth{n}}_{+}$.
\end{remark}

\begin{remark} \label{ideals-homogeneous} \rm The left ideals
$I_{\lambda_i}$, $i=0, \dots, l$, are graded by weight and by
$\lambda_j$-charge for $j=1, \dots, l$, and these gradings are
compatible.
\end{remark}

For any $\lambda \in P$ and character $\nu : Q \longrightarrow
\mathbb{C}^*$, we define a map $\tau_{\lambda,\nu}$ on
$\bar{\goth{n}}$ by
\[
\tau_{\lambda, \nu}(x_{\alpha}(m))=\nu(\alpha) x_{\alpha}(m-\langle
\lambda, \alpha \rangle)
\]
for $\alpha \in \Delta_{+}$ and $m \in \mathbb{Z}$.  It is easy to see
that $\tau_{\lambda,\nu}$ is an automorphism of $\bar{\goth{n}}$.
We will distinguish an important special case when $\nu$ is trivial
(i.e., $\nu=1$):  We set
\[
\tau_{\lambda}=\tau_{\lambda,1}.
\]

The map $\tau_{\lambda,\nu}$ extends canonically to an automorphism of
$U(\bar{\goth{n}})$, which we also denote by $\tau_{\lambda,\nu}$, so
that
\begin{equation}\label{def-lambda_i} 
\tau_{\lambda,\nu}(x_{\beta_1}(m_1) \cdots x_{\beta_k}(m_k))=
\nu(\beta_1+\cdots + \beta_k)
x_{\beta_1}(m_1-\langle \lambda, \beta_1 \rangle ) \cdots
x_{\beta_k}(m_k-\langle \lambda, \beta_k \rangle)
\end{equation} 
for $\beta_1, \dots, \beta_k \in \Delta_{+}$ and $m_1, \dots, m_k \in
\mathbb{Z}$.

Such automorphisms (or translations) generalize certain shift maps
that appeared in our previous work \cite{CalLM1}--\cite{CalLM2}
(recall \cite{CalLM1}, (3.20) and \cite{CalLM2}, (4.27)).  Notice that
for any $\lambda, \mu \in P$ and any characters $\nu$ and $\nu'$ on
$Q$, we have
\begin{equation}\label{tautau} 
\tau_{\lambda,\nu} \tau_{\mu,\nu'} = \tau_{\lambda+\mu,\nu \nu'} \; \;
\; \mbox{and} \; \; \;
\tau_{\lambda,\nu}^{-1}=\tau_{-\lambda,\nu^{-1}}.
\end{equation} 
Recall the (multiplicative) commutator map $c(\cdot, \cdot)$ on $P
\times P$, satisfying $c(\alpha,\beta) =(-1)^{\langle \alpha,\beta
\rangle}$ for $\alpha,\beta \in Q$.  For $\lambda \in P$, the map
\begin{equation}\label{c-aut} 
c_{\lambda}(\alpha)=c(\alpha,\lambda), \ \ \alpha \in Q,
\end{equation}
is a character on $Q$.
\begin{remark} \label{comparison}
\rm Assume that $a \in U(\bar{\goth{n}})$ is a nonzero element
homogeneous with respect to the weight and $\lambda_i$-charge
gradings. For any $\lambda \in P$ and character $\nu$ on $Q$,
$\tau_{\lambda, \nu}(a)$ and $\tau^{-1}_{\lambda, \nu}(a)$ are also
homogeneous and have the same $\lambda_i$-charge $n_i$ as $a$.
Moreover,
\begin{equation}\label{wttaulambdai}
\mbox{wt} \; \tau_{\lambda_i, \nu}(a) = \mbox{wt} \; a + n_{i},
\end{equation} 
\begin{equation}
\mbox{wt} \; \tau^{-1}_{\lambda_i, \nu}(a) = \mbox{wt} \; a -n_{i}.
\end{equation}
In particular, if $n_{i}>0$ then
\begin{equation}
\mbox{wt} \; \tau_{\lambda_i, \nu}(a) > \mbox{wt} \; a,
\end{equation} 
\begin{equation}
\mbox{wt} \; \tau^{-1}_{\lambda_i, \nu}(a) < \mbox{wt} \; a.
\end{equation}
\end{remark}

\begin{lemma} \label{lemma1}
For every $i = 1,\dots, l$ and character $\nu$ we have
\begin{equation} \label{inclusion_1}
\tau_{\lambda_i, \nu} (I_{\lambda_0}) \subset I_{\lambda_i}.
\end{equation}
\end{lemma}
\noindent {\it Proof:} Because $I_{\lambda_0}$ is a homogeneous ideal,
it is sufficient to consider $\tau_{\lambda_i}$. By (\ref{r}) and
(\ref{def-lambda_i}), for each $t \geq 2$ we have
\[
\tau_{\lambda_i}(R_{-1, t}^i)= R_{-1, t+2}^i+ a x_{\alpha_i}(-1),
\]
where $a \in U(\bar{\goth{n}})$, and 
\[
\tau_{\lambda_i} (R_{-1, t}^j) =R_{-1, t}^j \; \; \; \mbox{for} \; \;
j \neq i.
\]
Since $J$ is the left ideal of $U(\bar{\goth{n}})$ generated by the
$R_{-1, t}^i$ for $t \geq 2$ and $i=1, \dots, l$, we have
\[
\tau_{\lambda_i} (J) \subset I_{\lambda_i}.
\]
Note that for any $\beta \in \Delta_{+}$ and $m \geq 0$,
$x_{\beta}(m)$ can be expressed as a linear combination of monomials
of the form $x_{\alpha_{r_1}}(m_1) \cdots x_{\alpha_{r_k}}(m_k)$ such
that $m_k \geq 0$ and $\alpha_{r_1}, \dots, \alpha_{r_k}$ are simple
roots. If $r_k\neq i$ we have
\[
\tau_{\lambda_i}( x_{\alpha_{r_1}}(m_1) \cdots x_{\alpha_{r_k}}(m_k))
\in U(\bar{\goth{n}})\bar{\goth{n}}_{+}
\]
and if $r_k=i$ we have
\[
\tau_{\lambda_i}( x_{\alpha_{r_1}}(m_1) \cdots x_{\alpha_{r_k}}(m_k))
\in U(\bar{\goth{n}})\bar{\goth{n}}_{+}+
U(\bar{\goth{n}})x_{\alpha_i}(-1),
\]
whether $m_k>0$ or $m_k=0$. Thus
\[
\tau_{\lambda_i}(U(\bar{\goth{n}})\bar{\goth{n}}_{+}) \subset
I_{\lambda_i}.
\]
This completes the proof of (\ref{inclusion_1}). \; \; \; \; \; $\Box$

Consider the weights
\begin{equation}
\omega_i=\alpha_i-\lambda_i \in P
\end{equation}
and the automorphisms $\tau_{\omega_i}$ of $U(\bar{\goth{n}})$ for
$i=1, \dots, l$. These weights and certain maps associated with them
played an important role in \cite{Cal1} (see \cite{Cal1}, (4.12)).
For any character $\nu$ on $Q$, define the linear map
\begin{eqnarray}
\sigma_{\omega_i,\nu}: U(\bar{\goth{n}}) & \longrightarrow &
U(\bar{\goth{n}}) \nonumber \\ a & \mapsto &
\tau_{\omega_i,\nu}(a)x_{\alpha_i}(-1). \nonumber
\end{eqnarray}
This map is injective. To simplify the notation, we will write
\[
\sigma_{\omega_i} = \sigma_{\omega_i,1}.
\]

\begin{lemma} \label{lemma2}
For every $i=1, \dots, l$ and character $\nu$ we have
\begin{equation} \label{i-0}
\sigma_{\omega_i,\nu} (I_{\lambda_i}) \subset I_{\lambda_0}.
\end{equation}
\end{lemma}
\noindent {\it Proof:} Because the ideal $I_{\lambda_i}$ is generated
by homogenous elements, as in the previous lemma we may assume
that $\nu$ is trivial. Let $a \in U(\bar{\goth{n}})$. Then
\[
\sigma_{\omega_i}(ax_{\alpha_i}(-1))=\tau_{\omega_i}(a)x_{\alpha_i}(-2)
x_{\alpha_i}(-1)\in J \subset I_{\lambda_0}.
\]
Thus we have
\begin{equation} \label{i1}
\sigma_{\omega_i}(U(\bar{\goth{n}})x_{\alpha_i}(-1)) \subset
I_{\lambda_0}.
\end{equation}
Next we show that 
\begin{equation} \label{inclusion}
\sigma_{\omega_i}(x_{\beta}(m)) \in I_{\lambda_0}
\end{equation}
for any $\beta \in \Delta_{+}$ and $m \geq 0$. As in Lemma
\ref{lemma1} it is enough to show that (\ref{inclusion}) holds for
$x_{\alpha_j}(m)$, where $\alpha_j$ is a simple root and $m \geq 0$.
We shall use the fact that for any $\alpha, \beta \in \Delta$ such
that $\alpha+\beta \in \Delta$ we have $[x_{\alpha},
x_{\beta}]=C_{\alpha, \beta}x_{\alpha+\beta}$, where $C_{\alpha,
\beta} \neq 0$.  For any $m \in \mathbb{Z}$ we have
\begin{equation} \label{mainrel}
\sigma_{\omega_i}(x_{\alpha_j}(m))=\left\{ \begin{array}{cc}
x_{\alpha_i}(m-1)x_{\alpha_i}(-1) & {\rm if} \ i=j \\
x_{\alpha_j}(m+1)x_{\alpha_i}(-1) & {\rm if} \ i \neq j \; \;
\mbox{and} \; \; a_{ij}=-1 \\ x_{\alpha_j}(m)x_{\alpha_i}(-1) & {\rm
if} \ i \neq j \; \; \mbox{and} \; \; a_{ij}=0
\end{array} \right.
\end{equation}
For $m \geq 0$,
\begin{eqnarray}
&& x_{\alpha_i}(m-1)x_{\alpha_i}(-1)=x_{\alpha_i}(-1)x_{\alpha_i}(m-1)
\in I_{\lambda_0}, \nonumber \\ 
&&
x_{\alpha_j}(m+1)x_{\alpha_i}(-1)=
C_{\alpha_j,\alpha_i}x_{\alpha_i+\alpha_j}(m)+x_{\alpha_i}(-1)x_{\alpha_j}(m+1)
\in I_{\lambda_0} \; \; {\rm if} \ i \neq j, \ a_{ij}=-1, \nonumber \\
&& x_{\alpha_j}(m)x_{\alpha_i}(-1)=x_{\alpha_i}(-1)x_{\alpha_j}(m) \in
I_{\lambda_0} \; \; {\rm if} \ i \neq j, \ a_{ij}=0, \nonumber
\end{eqnarray}
and so $\sigma_{\omega_i}(x_{\alpha_j}(m)) \in I_{\lambda_0}$.  Thus
(\ref{inclusion}) holds and we have
\begin{equation} \label{i2}
\sigma_{\omega_i}(U(\bar{\goth{n}})\bar{\goth{n}}_{+}) \subset
I_{\lambda_0}.
\end{equation}

We now show that
\begin{equation} \label{maininclusion}
\sigma_{\omega_i}(J) \subset I_{\lambda_0}.
\end{equation}
Since $J$ is the left ideal of $U(\bar{\goth{n}})$ generated by
$R_{-1, t}^j$ for $t \in \mathbb{Z}$ and $j=1, \dots, l$, it is
sufficient to show that
\[
\sigma_{\omega_i}(R_{-1, t}^j) \in I_{\lambda_0} \; \; \mbox{for} \;
\; t \in \mathbb{Z}, \; j=1, \dots, l.
\]

For $i=j$, we have
\begin{equation} \label{uno}
\sigma_{\omega_i}(R_{-1, t}^i)= x_{\alpha_i}(-1)R_{-1,
t+2}^i+ax_{\alpha_i}(-1)^2 \in I_{\lambda_0},
\end{equation}
where $a \in U(\bar{\goth{n}})$.

If $i\neq j$ and $a_{ij}=0$, we have
\begin{equation} \label{duo}
\sigma_{\omega_i}(R_{-1, t}^j)=R_{-1,
t}^jx_{\alpha_i}(-1)=x_{\alpha_i}(-1)R_{-1, t}^j \in I_{\lambda_0}.
\end{equation}

If $i\neq j$ and $a_{ij}=-1$, then $\alpha_i+\alpha_j \in \Delta_{+}$
and $2\alpha_i+\alpha_j \notin \Delta$, and thus $[x_{\alpha_i}(m),
x_{\alpha_i+\alpha_j}(n)]=0$ for any $m,n \in \mathbb{Z}$.  Then for
any $t \in \mathbb{Z}$ we have
\begin{eqnarray} &&
\sigma_{\omega_i}(R_{-1, t}^{j})= \sum_{m_1,m_2 \leq-1, \;  m_1+m_2=-t}
x_{\alpha_j}(m_1+1)x_{\alpha_j}(m_2+1)x_{\alpha_i}(-1) \nonumber
\\
&& \label{sum1} = \sum_{m_1, m_2 \leq -1, \;  m_1+m_2=-t}
C_{\alpha_j,\alpha_i}
x_{\alpha_j}(m_1+1)x_{\alpha_i+\alpha_j}(m_2)
\\
&& \label{sum2} +
 \sum_{m_1, m_2 \leq -1, \;  m_1+m_2=-t}
C_{\alpha_j,\alpha_i}
x_{\alpha_i+\alpha_j}(m_1)x_{\alpha_j}(m_2+1) \\
&& \label{sum3} + x_{\alpha_i}(-1) \sum_{m_1, m_2 \leq -1, \; 
m_1+m_2=-t} x_{\alpha_j}(m_1+1)x_{\alpha_j}(m_2+1). 
\end{eqnarray}
The last term on the right-hand side is of the form
\begin{equation} \label{onee}
x_{\alpha_i}(-1)R_{-1, t-2}^j+a, \; \; \; \mbox{where} \; \; a \in
U(\bar{\goth{n}})\bar{\goth{n}}_{+}
\end{equation}
(and this is true even if $t=2$ or $t=3$) and thus this term is in
$I_{\lambda_0}$.  We rewrite (\ref{sum1}) as follows:
\begin{eqnarray} \label{two}
&& \sum_{m_1,  m_2
\leq -1, \; m_1+m_2=-t} C_{\alpha_j,\alpha_i}
x_{\alpha_j}(m_1+1)x_{\alpha_i+\alpha_j}(m_2) \nonumber \\ &&
\label{sum1-again} =C_{\alpha_j,\alpha_i}
x_{\alpha_j}(0)x_{\alpha_i+\alpha_j}(-t+1)+\sum_{m_1, m_2
\leq -1, \; m_1+m_2=-t+1} C_{\alpha_j,\alpha_i}
x_{\alpha_j}(m_1)x_{\alpha_i+\alpha_j}(m_2). 
\end{eqnarray}
Similarly,
\begin{eqnarray} \label{three}
&& \sum_{m_1, m_2 \leq
-1, \;  m_1+m_2=-t} C_{\alpha_j,\alpha_i}
x_{\alpha_i+\alpha_j}(m_1)x_{\alpha_j}(m_2+1) \nonumber \\
&& \label{sum2-again} =C_{\alpha_j,\alpha_i}
x_{\alpha_i+\alpha_j}(-t+1)x_{\alpha_j}(0)+ \sum_{m_1, m_2
\leq -1, \; m_1+m_2=-t+1}C_{\alpha_j,\alpha_i}
x_{\alpha_i+\alpha_j}(m_1)x_{\alpha_j}(m_2). 
\end{eqnarray}
Notice that
\begin{equation} \label{four}
[R_{-1, t-1}^j, x_{\alpha_i}(0)]
= R_{-1, t-1}^j x_{\alpha_i}(0) - x_{\alpha_i}(0) R_{-1, t-1}^j
\in J+U(\bar{\goth{n}})\bar{\goth{n}}_{+}=I_{\lambda_0}
\end{equation}
and that, on the other hand,
\begin{eqnarray} \label{five}
&& [R_{-1, t-1}^j, x_{\alpha_i}(0)] \\
&&=\sum_{m_1, m_2
\leq -1, \;  m_1+m_2=-t+1} C_{\alpha_j,\alpha_i}
x_{\alpha_j}(m_1)x_{\alpha_i+\alpha_j}(m_2) \nonumber \\
&&+ \sum_{m_1, m_2
\leq -1, \;  m_1+m_2=-t+1}C_{\alpha_j,\alpha_i}
x_{\alpha_i+\alpha_j}(m_1)x_{\alpha_j}(m_2).\nonumber
\end{eqnarray}
Combining (\ref{onee})--(\ref{five}) we obtain, for any 
$t \in \mathbb{Z}$ (in particular, for $t \geq 2$),
\begin{equation} \label{tre}
\sigma_{\omega_i}(R_{-1, t}^j) \in I_{\lambda_0} \; \; \; \mbox{if} \;
\; i \neq j \; \; \mbox{and} \; \; a_{ij}=-1.
\end{equation}

By (\ref{uno}), (\ref{duo}) and (\ref{tre}) we obtain
(\ref{maininclusion}). Combining (\ref{i1}), (\ref{i2}) and
(\ref{maininclusion}) we obtain the inclusion (\ref{i-0}). \; \; \; \;
\; $\Box$

Consider the composition $\sigma_{\omega_i}\tau_{\lambda_i}$.  As a
consequence of Lemmas \ref{lemma1} and \ref{lemma2} we obtain:

\begin{corollary} \label{one}
For every $i = 1,\dots, l$ and characters $\nu$ and $\nu'$ on $Q$ we
have
\begin{equation}
\sigma_{\omega_i,\nu} \tau_{\lambda_i, \nu'} (I_{\lambda_0}) 
\subset I_{\lambda_0}. \; \; \; \; \; \Box
\end{equation}
\end{corollary}

For any $\lambda \in P$ recall the linear isomorphism 
\[
e_{\lambda}: V_P \longrightarrow V_P.
\]
In particular, since
\[
e_{\lambda_i} \cdot v_{\lambda_0}=v_{\lambda_i}
\]
(recall (\ref{v_lambda_i}) and (\ref{epsilon-zero})), we have linear
isomorphisms
\begin{equation} \label{e_lambda_i}
e_{\lambda_i}:W(\lambda_0) \longrightarrow W(\lambda_i) \; \; \;
\mbox{for} \; \; \; i=1, \dots, l,
\end{equation}
given explictly as follows:  Since 
\[
e_{\lambda_i}x_{\alpha}(m)={c(\alpha, -\lambda_i)}
x_{\alpha}(m-\langle \alpha, \lambda_i \rangle)e_{\lambda_i} \; \; \;
\mbox{for} \; \; \alpha \in \Delta_{+} \; \; \mbox{and} \; \; m
\in \mathbb{Z}
\] 
(recall (\ref{almost-comm})), we have
\begin{equation}\label{etau}
e_{\lambda_i}(a \cdot v_{\lambda_0})=\tau_{\lambda_i,c_{-\lambda_i}}(a) \cdot
v_{\lambda_i}, \; \; \; a \in U(\bar{\goth{n}}).
\end{equation}
Generalizing the corresponding constructions in
\cite{CalLM1}--\cite{CalLM2}, we construct a linear lifting
\begin{equation}
\widehat{e_{\lambda_i}}:W^N(\lambda_0) \longrightarrow W^N(\lambda_i)
\end{equation}
of each map (\ref{e_lambda_i}),
in the sense that the diagram
\[ \CD W^N(\lambda_0) @> \widehat{e_{\lambda_i}}> >
W^N(\lambda_i) \\ @V\pi_{\lambda_{0}}VV
@V\pi_{\lambda_{i}}VV \\ W(\lambda_0)\ @> {e_{\lambda_i}}> >
W(\lambda_i)  \endCD
\]
will commute, for $i=1, \dots, l$, by taking
\begin{equation} \label{formula}
\widehat{e_{\lambda_i}}(a \cdot v_{\lambda_0}^N)=\tau_{\lambda_i,c_{-\lambda_i}}(a)
\cdot v^N_{\lambda_i} \; \; \; \mbox{for} \; \; a \in U(\bar{\goth{n}}_{-});
\end{equation}
the map $\widehat{e_{\lambda_i}}$ is well defined since
$W^N(\lambda_0)$ is a free $U(\bar{\goth{n}}_{-})$-module.

We have:

\begin{corollary} \label{corollary_two}
For $i = 1,\dots, l$,
\begin{equation} \label{hat-inclusion}
\widehat{e_{\lambda_i}}(I_{\lambda_0} \cdot v^N_{\lambda_0}) \subset
I_{\lambda_i} \cdot v_{\lambda_i}^N.
\end{equation}
\end{corollary}
\noindent{\it Proof:} Let $a \in I_{\lambda_0}$. We may assume that $a
\in U(\bar{\goth{n}}_{-})$. Indeed, since $a \in U(\bar{\goth{n}})$ we
can write $a=b+c$, where $b \in U(\bar{\goth{n}})\bar{\goth{n}}_{+}$
and $c \in U(\bar{\goth{n}}_{-})$ (recall (\ref{decomp})), and since
$b \in I_{\lambda_0}$, $c \in I_{\lambda_0}$ as well.  Now by applying
(\ref{formula}) and Lemma \ref{lemma1} to the element $a \in
I_{\lambda_0} \cap U(\bar{\goth{n}}_{-})$ we obtain the inclusion
(\ref{hat-inclusion}).  \; \; \; \; \; $\Box$

Recall the weights $\omega_i=\alpha_i-\lambda_i$, $i=1, \dots, l$. We
now consider the linear isomorphism
\begin{equation}
e_{\omega_i}: V_P\longrightarrow V_P
\end{equation}
and its restriction to the principal-like subspace $W(\lambda_i)$.
Since, by (\ref{almost-comm}),
\[
e_{\omega_i}x_{\alpha}(m)={c(\alpha, -\omega_i)}
x_{\alpha}(m-\langle \alpha, \omega_i \rangle)e_{\omega_i} \; \; \;
\mbox{for} \; \; \alpha \in \Delta_{+} \; \; \mbox{and} \; \; m
\in \mathbb{Z}
\] 
and
\[
e_{\omega_i} \cdot v_{\lambda_i} = e_{\omega_i} \cdot e^{\lambda_i} =
\epsilon(\omega_i,\lambda_i)e^{\alpha_i} =
\epsilon(\omega_i,\lambda_i)x_{\alpha_i}(-1) \cdot v_{\lambda_0},
\] 
we have
\begin{equation}
e_{\omega_i}(a \cdot v_{\lambda_i})= \epsilon(\omega_i,\lambda_i)
\tau_{\omega_i, c_{-\omega_i}}(a)x_{\alpha_i}(-1) \cdot v_{\lambda_0}
= \epsilon(\omega_i,\lambda_i) \sigma_{\omega_i, c_{-\omega_i}}(a)
\cdot v_{\lambda_0}
\end{equation}
for $a \in U(\bar{\goth{n}})$, giving us a linear injection 
\begin{equation}  \label{e_omega_i}
e_{\omega_i}:W(\lambda_i) \longrightarrow W(\lambda_0)
\end{equation}
for $i=1, \dots, l$.

\section{Presentations of the principal subspaces}
\setcounter{equation}{0}

In this section we prove natural presentations of the principal
subspaces of the level one standard modules for $\widehat{\goth{g}}$,
where $\goth{g}$ is of type $A_l$, $l \geq 1$, $D_l$, $l \geq 4$,
$E_6$, $E_7$ and $E_8$ (Theorem \ref{presentations-ps}).  More
generally, we prove natural presentations of the principal-like
subspaces $W(\lambda_i) \subset V_P$ for $i=0, \dots, l$ (Theorem
\ref{presentations-finite}).  In order to obtain these results we
generalize the proof of Theorem 2.2 (which is equivalent to Theorem
2.1) of \cite{CalLM1}. Following \cite{Cal1}, we give a reformulation
of the presentations of the principal subspaces, and more generally,
of the principal-like subspaces, in terms of ideals of a certain
completion of the universal enveloping algebra $U(\bar{\goth{n}})$
(Theorem \ref{presentations-infinite}), and we show that Theorem
\ref{presentations-finite} implies Theorem
\ref{presentations-infinite}.

First we state:

\begin{theorem} \label{presentations-finite}
For every $i = 0,\dots, l$, we have
\begin{equation}
\mbox{\rm Ker} \; f_{\lambda_i}=I_{\lambda_i},
\end{equation}
or equivalently,
\begin{equation} \label{main_formula}
\mbox{\rm Ker} \; \pi_{\lambda_i} = I_{\lambda_i} \cdot v_{\lambda_i}^N.
\end{equation}
\end{theorem}

(To see that (\ref{main_formula}) implies that $\mbox{\rm Ker} \;
f_{\lambda_i} \subset I_{\lambda_i}$, we use the direct sum
decomposition (\ref{decomp}).)

Recall that when $\langle \Lambda_i, {\bf k} \rangle =1$, the
principal-like subspace $W(\lambda_i)$ agrees with the principal
subspace $W(\Lambda_i)$ of the level one standard module
$L(\Lambda_i)$, and $f_{\lambda_i}$ and $\pi_{\lambda_i}$ agree with
the maps $f_{\Lambda_i}$ and $\pi_{\Lambda_i}$, respectively
(cf. Remarks \ref{w=w} and \ref{f=f}).  As a particular case of
Theorem \ref{presentations-finite} we have the following presentations
of the principal subspaces $W(\Lambda_i)$:

\begin{theorem} \label{presentations-ps}
For any fundamental weight $\Lambda_i$ such that $\langle \Lambda_i,
{\bf k} \rangle =1$,
\begin{equation}
\mbox{ \rm Ker} \; f_{\Lambda_i}=I_{\lambda_i},
\end{equation}
or equivalently,
\begin{equation} 
\mbox{\rm Ker} \; \pi_{\Lambda_i} = I_{\lambda_i} \cdot
v_{\Lambda_i}^N.
\end{equation}
\end{theorem}

\noindent {\it Proof of Theorem \ref{presentations-finite}:} For each
$j = 1,\dots, l$, the square of the vertex operator $Y(e^{\alpha_j},
x)$ is well defined and equals zero on $V_P$ and hence on each
$W(\lambda_i)$.  Thus $Y(e^{\alpha_j}, x)^2=\sum_{t \in
\mathbb{Z}}R_t^jx^{t-2}$ equals zero on $W(\lambda_i)$ (recall
(\ref{R})), so that $R_t^j = 0$ on $W(\lambda_i)$ for each $t \in
\mathbb{Z}$.  This combined with the highest weight vector property of
$v_{\lambda_i}$ and the fact that
\[
x_{\alpha_i}(-1) \cdot v_{\lambda_i} = 0 \; \; \; \mbox{if} \; \; i>0
\]
(which follows from (\ref{vertex})) implies that
\[I_{\lambda_i} \subset \mbox{Ker} \; f_{\lambda_i},
\]
and so
\[
I_{\lambda_i} \cdot v_{\lambda_i}^N \subset \mbox{Ker} \;
\pi_{\lambda_i}
\]
for $i=0, \dots, l$. We will now prove the inclusions
\begin{equation} \label{inc}
\mbox{Ker} \; \pi_{\lambda_i} \subset I_{\lambda_i} \cdot
v_{\lambda_i}^N \; \; \; \mbox{for} \; \; i=0, \dots, l.
\end{equation}

We first claim that (\ref{main_formula}) with $i=1, \dots, l$ follows
{}from (\ref{main_formula}) with $i=0$, whose truth we now assume.
Let $i > 0$ and let $a \cdot v_{\lambda_i}^N \in \mbox{Ker} \;
\pi_{\lambda_i}$, where $a \in U(\bar{\goth{n}})$. Then $a \cdot
v_{\lambda_i}=0$ in $W(\lambda_i)$, and so
\[
\tau^{-1}_{\lambda_i,c_{-\lambda_i}} (a) \cdot v_{\lambda_0}=0
\]
in $W(\lambda_0)$, by (\ref{etau}).  Thus
\[
\tau^{-1}_{\lambda_i, c_{-\lambda_i}}(a) \cdot v_{\lambda_0}^N \in
\mbox{Ker} \; \pi_{\lambda_0}=I_{\lambda_0} \cdot v_{\lambda_0}^N.
\]
Now by Corollary \ref{corollary_two} we have
\[
\widehat{e_{\lambda_i}} (\tau^{-1}_{\lambda_i, c_{-\lambda_i}}(a) \cdot
v_{\lambda_0}^N) \in I_{\lambda_i} \cdot v_{\lambda_i}^N,
\]
and from (\ref{formula}) we get
\[
a \cdot v_{\lambda_i}^N \in I_{\lambda_i} \cdot v_{\lambda_i}^N \; \;
\; \mbox{for} \; \; \; i=1, \dots , l.
\]
This proves our claim.

As in \cite{CalLM1}--\cite{CalLM2} we will use a contradiction
argument to prove the inclusion
\begin{equation} \label{inc_zero}
\mbox{Ker} \; \pi_{\lambda_0}\subset I_{\lambda_0} \cdot
v_{\lambda_0}^N,
\end{equation}
which is all that remains to prove.  Suppose then that there exists $a
\in U(\bar{\goth{n}})$ such that
\begin{equation} \label{contradiction}
a \cdot v_{\lambda_0}^N \in \mbox{Ker} \; \pi_{\lambda_0} \;\; \;
\mbox{but} \; \; \; a \cdot v_{\lambda_0}^N \notin I_{\lambda_0} \cdot
v_{\lambda_0}^N.
\end{equation}
We may and do assume that $a$ is homogeneous with respect to the
weight and $\lambda_i$-charge gradings (recall Remarks \ref{L(0)} and
\ref{ideals-homogeneous}).  The element $a$ is certainly nonzero, and
it is also nonconstant because otherwise, $a \cdot v_{\lambda_0}^N
\notin \mbox{Ker} \; \pi_{\lambda_0}$.  Using this and the
decomposition (\ref{decomp}), we further see that $a$ has positive
weight, since otherwise, $a$ would be an element of
$U(\bar{\goth{n}})\bar{\goth{n}}_{+}$ and hence of $I_{\lambda_0}$.
Denote by $n$ the total charge of $a$, namely, $n=n_1+ \cdots +n_l$,
where $n_i \geq 0$ is the $\lambda_i$-charge of $a$ for $i=1, \dots,
l$. Note that $n>0$; otherwise, $a$ is constant.  Take $n$ to be the
minimum total charge for all the homogeneous elements $a$ satisfying
(\ref{contradiction}).  Among all the homogeneous elements of total
charge $n$ satisfying (\ref{contradiction}), we choose $a$ to be an
element of the smallest possible (necessarily positive) weight.  Fix
any index $i$ for which $n_i > 0$.

We claim that $a \in I_{\lambda_i}$.  Assume then that
\begin{equation} \label{noclaim1}
a \notin I_{\lambda_i}.
\end{equation}
Since $a \cdot v_{\lambda_0}^N \in \mbox{Ker} \; \pi_{\lambda_0}$, $a
\in \mbox{Ker} \; f_{\lambda_0} \subset \mbox{Ker} \; f_{\lambda_i}$
(see (\ref{inclusion_kernels})). Then $a \cdot v_{\lambda_i} =0$ in
$W(\lambda_i)$, and so
\[
e_{\lambda_i} (\tau^{-1}_{\lambda_i, c_{-\lambda_i}} (a) \cdot
v_{\lambda_0})=a \cdot v_{\lambda_i}=0,
\]
by (\ref{etau}).  This implies that
$\tau^{-1}_{\lambda_i,c_{-\lambda_i}}(a) \cdot v_{\lambda_0}=0$, or
equivalently, that
\begin{equation} \label{cont1} \tau^{-1}_{\lambda_i,
c_{-\lambda_i}}(a) \cdot v_{\lambda_0}^N \in \mbox{Ker} \;
\pi_{\lambda_0}.
\end{equation}
Also,
\begin{equation} \label{cont2} \tau^{-1}_{\lambda_i,
c_{-\lambda_i}}(a) \cdot v_{\lambda_0}^N \notin I_{\lambda_0}
\cdot v_{\lambda_0}^N;
\end{equation}
otherwise, we would have $\tau^{-1}_{\lambda_i, c_{-\lambda_i}}(a) \in
I_{\lambda_0}$ (from (\ref{decomp})) and so by Lemma \ref{lemma1} we
would get $a \in I_{\lambda_i}$, contradicting (\ref{noclaim1}).  The
elements $\tau^{-1}_{\lambda_i, c_{-\lambda_i}}(a)$ and $a$ have the
same $\lambda_j$-charge for each $j$, and hence the same total charge
$n$, and
\begin{equation} \label{cont3} \mbox{wt} \; \tau^{-1}_{\lambda_i,
c_{-\lambda_i}}(a) < \mbox{wt} \; a,
\end{equation}
by Remark \ref{comparison}.  Now (\ref{cont1}), (\ref{cont2}) and
(\ref{cont3}) contradict our choice of $a$ and thus we get
\begin{equation}
a \in I_{\lambda_i}.
\end{equation}

Since $I_{\lambda_i}=I_{\lambda_0}+U(\bar{\goth{n}})x_{\alpha_i}(-1)$,
there exist $b \in I_{\lambda_0}$ and $c \in U(\bar{\goth{n}})$ such
that
\begin{equation} \label{decomposition}
a=b+cx_{\alpha_i}(-1).
\end{equation}
We may and do assume that $b$ and $c$ are homogeneous with respect to
the weight and $\lambda_j$-charge gradings.  Note that $b$ has the
same weight and the same total charge $n$ as $a$; the total charge of
$c$ is $n-1$; and also, $\mbox{wt} \; c = \mbox{wt} \; a - 1$.

We now claim that
\begin{equation} \label{claim2}
cx_{\alpha_i}(-1) \in I_{\lambda_0}.
\end{equation}
Assume that (\ref{claim2}) does not hold. Then 
\begin{equation} \label{contra1} \tau^{-1}_{\alpha_i,
c_{-\alpha_i}}(c) \cdot v_{\lambda_0}^N \notin I_{\lambda_0} \cdot
v_{\lambda_0}^N.
\end{equation}
Indeed, from (\ref{tautau}) we have
\[
\tau_{\lambda_i, c_{-\alpha_i}}\tau^{-1}_{\alpha_i, c_{-\alpha_i}}
=\tau_{\lambda_i, c_{-\alpha_i}}\tau_{-\alpha_i, (c_{-\alpha_i})^{-1}}
=\tau_{-\omega_i}
\]
and
\[
\sigma_{\omega_i}\tau_{-\omega_i}(c)
=\tau_{\omega_i}(\tau_{-\omega_i}(c))x_{\alpha_i}(-1)
=cx_{\alpha_i}(-1),
\]
so that
\[
cx_{\alpha_i}(-1)= (\sigma_{\omega_i} \tau_{\lambda_i, c_{-\alpha_i}})
(\tau^{-1}_{\alpha_i, c_{-\alpha_i}}(c)),
\]
and Corollary \ref{one} yields $\tau^{-1}_{\alpha_i, c_{-\alpha_i}}(c)
\notin I_{\lambda_0}$; thus (\ref{contra1}) holds.  By
(\ref{almost-comm}), for $f \in U(\bar{\goth{n}})$,
\[
e_{\alpha_i}f = \tau_{\alpha_i, c_{-\alpha_i}}(f)e_{\alpha_i}
\]
as operators, so that
\[
e_{\alpha_i} (\tau^{-1}_{\alpha_i, c_{-\alpha_i}}(c) \cdot
v_{\lambda_0}) =c e_{\alpha_i} \cdot v_{\lambda_0} =c e^{\alpha_i} =c
x_{\alpha_i}(-1)\cdot v_{\lambda_0}= (a-b) \cdot v_{\lambda_0}=0,
\]
and so
\begin{equation} \label{contra2}
\tau^{-1}_{\alpha_i, c_{-\alpha_i}}(c) \cdot v_{\lambda_0}^N \in \mbox{Ker} \;
\pi_{\lambda_0}.
\end{equation}
Also, by Remark \ref{comparison}, $\tau^{-1}_{\alpha_i,
c_{-\alpha_i}}(c)$ is homogeneous and has the same total charge as
$c$, namely, $n-1$.  Thus (\ref{contra1}) and (\ref{contra2})
contradict our choice of the element $a$, and so we have
$cx_{\alpha_i}(-1) \in I_{\lambda_0}$, proving our claim
(\ref{claim2}).

It follows that $a \in I_{\lambda_0}$.  This shows that our initial
assumption is false, and therefore we have
\[
\mbox{Ker} \; \pi_{\lambda_0} \subset I_{\lambda_0} \cdot
v_{\lambda_0}^N. \; \; \; \; \; \Box
\]

\begin{remark} \em We now comment on the similarities between the
proofs of Theorem \ref{presentations-finite} in this paper and Theorem
2.2 of \cite{CalLM1}, which deals with the $\goth{sl}(2)$ level one
case.  Note that when $\goth{g}=\goth{sl}(2)$ the principal-like
subspaces coincide with the principal subspaces.  One difference
between these two proofs is that here we use $U(\bar{\goth{n}})$
rather than $U(\bar{\goth{n}}_{-})$ (which was sufficient for the
$\goth{sl}(2)$ case because $\bar{\goth{n}}$ is abelian), and this
simplifies the argument somewhat.  In the proof of Theorem
\ref{presentations-finite} the claim that $a \in I_{\lambda_i}$ shows
that there is in fact a homogeneous element lying in
$U(\bar{\goth{n}})x_{\alpha_i}(-1)$ and in addition having all the
properties of $a$, since $cx_{\alpha_i}(-1)$ has the same weight and
charge as the element $a$ and, in addition, it satisfies
(\ref{contradiction}).  When $\goth{g}=\goth{sl}(2)$ this claim is
similar to the one that appears in the proof of Theorem 2.2 of
\cite{CalLM1}, with $U(\bar{\goth{n}})$ replaced by
$U(\bar{\goth{n}}_{-})$.  Also, the proof of (\ref{claim2}) follows
the lines of the last part of the proof of Theorem 2.2 of
\cite{CalLM1}, except that here our minimal counterexample involves
charge as well as weight.  In the $\goth{sl}(2)$ case,
$\tau^{-1}_{\alpha_1, c_{-\alpha_1}}$ and $e_{\alpha_1}$,
respectively, in this paper are the same as $\tau^{-2}$ and
$e^{\alpha/2} \circ e^{\alpha/2}$, respectively, in \cite{CalLM1}, and
so in the $\goth{sl}(2)$ case, it was not necessary to minimize charge
as well as weight.
\end{remark}

Recall from \cite{Cal1} the two-sided ideal, denoted by ${\mathcal
J}$, of $\widetilde{U(\bar{\goth{n}})}$, the completion of
$U(\bar{\goth{n}})$ in the sense of \cite{LW3} or \cite{MP1},
generated by the elements $R_{t}^{j}$ for all $t \in \mathbb{Z}$ and
$j=1, \dots, l$, where the $R_t^j$ are defined in (\ref{R}).
The decomposition (\ref{decomp}) implies:
\begin{equation} \label{decomptilde}
\widetilde{U(\bar{\goth{n}})} = U(\bar{\goth{n}}_{-}) +
\widetilde{U(\bar{\goth{n}})\bar{\goth{n}}_{+}}.
\end{equation}

We now give a different description of the ideals $\mbox{Ker} \;
f_{\Lambda_i}$, and more generally, $\mbox{Ker} \; f_{\lambda_i}$, in
terms of the two-sided ideal ${\mathcal J}$:

\begin{theorem}  \label{presentations-infinite}
The annihilator $\mbox{\rm Ker} \; f_{\Lambda_0}$ in
$U(\bar{\goth{n}})$ of the highest weight vector of $L(\Lambda_0)$ is
described as follows:
\begin{equation} \label{equation_1}
\mbox{\rm Ker} \; f_{\Lambda_0} \equiv {\mathcal J} \; \; \;
\mbox{\rm modulo} \; \; \;
\widetilde{U(\bar{\goth{n}})\bar{\goth{n}}_{+}}.
\end{equation}
Moreover, for the principal-like subspaces $W(\lambda_i)$, $i =
1,\dots, l$, we have:
\begin{equation} \label{equation_2}
\mbox{\rm Ker} \; f_{\lambda_i}\equiv {\mathcal J} +
U(\bar{\goth{n}})x_{\alpha_i}(-1) \; \; \; \mbox{\rm modulo} \; \; \;
\widetilde{U(\bar{\goth{n}})\bar{\goth{n}}_{+}}.
\end{equation}
\end{theorem}
\noindent {\em Proof.} The inclusions
\[
{\mathcal J}  \subset \mbox{Ker} \; f_{\Lambda_0}
\; \; \; \mbox{modulo} \; \; \; 
\widetilde{U(\bar{\goth{n}})\bar{\goth{n}}_{+}}
\]
and
\[
{\mathcal J} + U(\bar{\goth{n}})x_{\alpha_i}(-1) \subset \mbox{Ker} \;
f_{\lambda_i} \; \; \; \mbox{modulo} \; \; \;
\widetilde{U(\bar{\goth{n}})\bar{\goth{n}}_{+}}
\]
follow from the definition of the ideal ${\mathcal J}$, the first
paragraph of the proof of Theorem \ref{presentations-finite}, and
(\ref{decomptilde}).

Observe that
\begin{equation} \label{small}
J \subset {\mathcal J}  \; \; \; \mbox{modulo} \; \; \; 
\widetilde{U(\bar{\goth{n}})\bar{\goth{n}}_{+}}.
\end{equation}

Let $a \in \mbox{Ker} \; f_{\Lambda_0}$. By Theorem
\ref{presentations-ps} and (\ref{I_0}) we have $a \in J
+U(\bar{\goth{n}})\bar{\goth{n}}_{+}$. Then by (\ref{small}) we have
$a \in {\mathcal J} \; \; \; \mbox{modulo} \; \; \;
\widetilde{U(\bar{\goth{n}})\bar{\goth{n}}_{+}}$. This proves the
inclusion
\[
\mbox{\rm Ker} \; f_{\Lambda_0} \subset {\mathcal J} \; \; \;
\mbox{modulo} \; \; \;
\widetilde{U(\bar{\goth{n}})\bar{\goth{n}}_{+}},
\]
and thus (\ref{equation_1}).  Now let $a \in \mbox{Ker} \;
f_{\lambda_i}$, $i \neq 0$.  Using Theorem \ref{presentations-finite},
(\ref{I_lambda_i}) and (\ref{small}) we obtain the inclusion
\[
\mbox{\rm Ker} \; f_{\lambda_i} \subset {\mathcal J} +
U(\bar{\goth{n}})x_{\alpha_i}(-1) \; \; \; \mbox{modulo} \; \; \;
\widetilde{U(\bar{\goth{n}})\bar{\goth{n}}_{+}},
\]
and so formula (\ref{equation_2}) holds. \; \; \; \; \; $\Box$

\section{$q$-difference equations}
\setcounter{equation}{0}

In this section we will use the presentations proved in Theorems
\ref{presentations-finite} and \ref{presentations-ps} in order to
construct canonical exact sequences for the principal(-like)
subspaces.  As a consequence we derive $q$-difference equations
(recursions) for the (multi-)graded dimensions of these subspaces and
by solving these equations we obtain explicit formulas for the graded
dimensions of the subspaces.  Since the proofs of the results in this
section use exactly the same arguments as the proofs of the
corresponding results in \cite{CLM1}, \cite{CLM2} and \cite{Cal1}, we
will only sketch the proofs here.

Using the injection (\ref{e_omega_i}) and the surjection
(\ref{intertwining}), associated to the relevant intertwining
operators, we obtain:

\begin{theorem} \label{1st}  For every $i = 1,\dots, l$
we have the following short exact sequence of maps among
principal-like subspaces:
\begin{equation} \label{ses1} \CD 0 @> > >W(\lambda_i) @> e_{\omega_i} > >
W(\lambda_0) @> \mathcal{Y}_c(e^{\lambda_i},x) >> W(\lambda_i) @> >> 0 .
\endCD
\end{equation}
In particular, for each fundamental weight $\Lambda_i$ of
$\widehat{\goth{g}}$ such that $\langle \Lambda_i, {\bf k} \rangle=1$
we have the following short exact sequence among principal subspaces:
\begin{equation} 
\CD 0 @> > >W(\Lambda_i) @> e_{\omega_i} > >
W(\Lambda_0) @> \mathcal{Y}_c(e^{\lambda_i},x) >> W(\Lambda_i) @> >> 0 .
\endCD
\end{equation}
\end{theorem}
\noindent {\em Proof:} The proof follows the same lines as in
\cite{CLM1}--\cite{CLM2}, or more specifically, as in \cite{Cal1}.
The main step is the exactness of the chain complex at the middle, and
this is where the presentation result, Theorem
\ref{presentations-finite} (or Theorem \ref{presentations-ps}), is
used.  In fact, the key steps, easily proved from the above, are:
\[
\mbox{\rm Ker} \; {\mathcal Y}_c(e^{\lambda_i}, x) =
\mbox{\rm Ker} \; f_{\lambda_i} \cdot v_{\lambda_0},
\]
so that by Theorem \ref{presentations-finite},
\[
\mbox{\rm Ker} \; {\mathcal Y}_c(e^{\lambda_i}, x) =
I_{\lambda_i} \cdot v_{\lambda_0},
\]
whereas
\[
\mbox{\rm Im} \; e_{\omega_i} =U(\bar{\goth{n}})x_{\alpha_i}(-1) \cdot
v_{\lambda_0},
\]
and these spaces agree, by (\ref{I_lambda_i}).  \; \; \; \; \; $\Box$

Combining the short exact sequence (\ref{ses1}) with two additional
maps, which are isomorphisms (recall (\ref{e_lambda_i})), we obtain
the diagram
\[\xymatrix{ 0 \ar[r] & W(\Lambda_0) \ar^{e_{\lambda_i}}[d] \ar@{.>}[rd] & & &  \\
0 \ar[r] & W(\lambda_i) \ar^{e_{\omega_i}}[r] & W(\Lambda_0)
\ar^{\mathcal{Y}_c(e^{\lambda_i},x)}[r] \ar@{.>}[rd] & W(\lambda_i)
\ar^{e_{-\lambda_i}}[d] \ar[r] & 0 \\ & & & W(\Lambda_0) \ar[r] & 0,
}\]
and thus another formulation of Theorem \ref{1st}, which uses only the
principal subspace $W(\Lambda_0)$ (note that $e_{\alpha_i}$ is
proportional to $e_{\omega_i}e_{\lambda_i}$, by
(\ref{e-multiplication})):

\begin{theorem} \label{2nd}
For every $i = 1,\dots, l$ we have the following short exact sequence:
\[ \CD 0 @> > >W(\Lambda_0) @> e_{\alpha_i} > >
W(\Lambda_0) @> e_{-\lambda_i} \circ \mathcal{Y}_c(e^{\lambda_i},x) >>
W(\Lambda_0) @> >> 0 .  \endCD
\; \; \; \; \; \Box \]
\end{theorem}

As we recalled in Section 2, the vector space $V_P$ and its subspaces
$W(\Lambda_i)$ and $W(\lambda_i)$ are graded by weight and charge, and
these gradings are compatible. We consider the multi-graded dimension
of the principal-like subspace $W(\lambda_i)$, for each $i=0,\dots,
l$:
\[\chi_{W(\lambda_i)}(x_1,\dots,x_l;q)={\rm tr}_{W(\lambda_i)}
x_1^{\lambda_1} \cdots x_l^{\lambda_l} q^{L(0)},
\]
where $x_1, \dots, x_l$ and $q$ are commuting formal variables and
$L(0)$ is the standard Virasoro algebra operator, introduced earlier.
Note that
\[
\chi_{W(\Lambda_0)}(x_1, \dots, x_l;q) \in \mathbb{C}[[x_1, \dots, x_l,q]].
\]
As in \cite{CLM2}, in order to avoid the multiplicative factors
$x_1^{\langle \lambda_1, \lambda_i \rangle} \cdots x_l^{ \langle
\lambda_l, \lambda_i \rangle}q^{\frac{1}{2} \langle \lambda_i,
\lambda_i \rangle}$, we use the following slightly modified graded
dimensions:
\begin{equation}\label{chi'chi}
\chi'_{W(\lambda_i)}(x_1,\dots,x_l;q) =x_1^{-\langle
\lambda_i,\lambda_1 \rangle} \cdots x_l^{-\langle \lambda_i,\lambda_l
\rangle} q^{-\langle \lambda_i, \lambda_i \rangle/2}
\chi_{W(\lambda_i)}(x_1,\dots,x_l;q)
\end{equation}
for $i = 1,\dots, l$, so that $\chi_{W(\lambda_i)}'(x_1, \dots, x_l;q)
\in \mathbb{C}[[x_1, \dots, x_l,q]]$ (and
$\chi_{W(\lambda_0)}'=\chi_{W(\Lambda_0)}$).

Using the isomorphism (\ref{e_lambda_i})--(\ref{etau}) and
(\ref{wttaulambdai}), we easily obtain
\begin{equation}\label{chi'}
\chi'_{W(\lambda_i)}(x_1, \dots, x_l;q)=
\chi_{W(\Lambda_0)}(x_1, \dots, x_iq, \dots, x_l;q)
\end{equation}
for $i>0$.

Following \cite{Cal1} (see also \cite{CLM1} and \cite{CLM2}), from
Theorem \ref{1st} (or equivalently from Theorem \ref{2nd}) we obtain a
canonical system of $q$-difference equations for the graded dimension
of the principal subspace $W(\Lambda_0)$:

\begin{theorem} \label{q-diff-eq} We have the
following system of $q$-difference equations for $i = 1,\dots, l$:
\[
\chi_{W(\Lambda_0)}(x_1,\dots,x_l;q)=\chi_{W(\Lambda_0)}(x_1,\dots,x_iq,\dots,x_l;q)+
(x_i q) \chi_{W(\Lambda_0)}(x_1q^{m_{i1}},\dots,x_lq^{m_{il}};q),
\]
where $M=(m_{ij})_{1 \leq i, j \leq
l}$ is the Cartan matrix of $\goth{g}$.
\; \; \; \; \; $\Box$
\end{theorem}

\begin{remark}
\em The first term on the right-hand side comes from the map
$\mathcal{Y}_c(e^{\lambda_i},x)$ combined with (\ref{chi'}); the
factor $x_i q$ in the second term comes from the fact that
$e_{\alpha_i} \cdot v_{\lambda_0} =e^{\alpha_i}$, which has
$\lambda_j$-charge $\delta_{ij}$ and weight $1$; and the expressions
$x_iq^{m_{ij}}$ come from (\ref{def-lambda_i}) with
$\lambda=\alpha_i$, together with the analogue
\[
e_{\alpha_i}(a \cdot v_{\lambda_0})=\tau_{\alpha_i,c_{-\alpha_i}}(a) \cdot
e^{\alpha_i}, \; \; \; a \in U(\bar{\goth{n}})
\]
of (\ref{etau}).
\end{remark}

For any $n_1, \dots, n_l \geq 0$ we define $f_{n_1, \dots, n_l}(q) \in
\mathbb{C}[[q]]$ by:
\[
\chi_{W(\Lambda_0)}(x_1,\dots,x_n,q)=\sum_{n_1,\dots, n_l \geq 0}
f_{n_1,\dots,n_l}(q) x_1^{n_1} \cdots x_l^{n_l}.
\]
As in \cite{Cal1}, it is straightforward to show that the system of
$q$-difference equations obtained in Theorem \ref{q-diff-eq} has a
unique solution in $\mathbb{C}[[x_1,\dots,x_l,q]]$ with the initial
condition $f_{0,\dots,0}(q)=1$ ($f_{0, \dots, 0}(q)$ being the graded
dimension of the subspace consisting of the elements of charge zero
with respect to $\lambda_1, \dots, \lambda_l$), and to generate the
solution.  As usual, for any nonnegative integer $n$ set
\[
(q)_n=(1-q) \cdots (1-q^n).
\]

Following the same argument as in the proof of Corollary 4.1 in
\cite{Cal1}, we obtain from Theorem \ref{q-diff-eq}, using the $l$
equations in succession, the (multi-)graded dimension of the principal
subspace $W(\Lambda_0)$; then we invoke (\ref{chi'}) to obtain the
(multi-)graded dimensions of the principal-like subspaces
$W(\lambda_i)$ for $i > 0$:

\begin{corollary} \label{characters} Let $M=(m_{ij})_{1 \leq i, j \leq
l}$ be the Cartan matrix of $\goth{g}$.  We have
\begin{equation} \label{char-0}
\chi_{W(\Lambda_0)}(x_1,\dots,x_l;q)=\sum_{{\bf n}=(n_1,\dots,n_l)}
\frac{q^{\frac{{\bf n} M {\bf n}^T}{2}}}{(q)_{n_1} \cdots
(q)_{n_l}}x_1^{n_1} \cdots x_l^{n_l}
\end{equation}
(where each $n_j \geq 0$).  Moreover, for $i = 1,\dots, l$,
\begin{equation}\label{char-i} 
\chi'_{W(\lambda_i)}(x_1,\dots,x_l;q)=\sum_{{\bf
n}=(n_1,\dots,n_l)} \frac{q^{\frac{{\bf n} M {\bf
n}^T+2n_i}{2}}}{(q)_{n_1} \cdots (q)_{n_l}}x_1^{n_1} \cdots x_l^{n_l}.
\; \; \; \; \; \Box
\end{equation}
\end{corollary}

\begin{remark}
\em Combining (\ref{char-i}) with (\ref{chi'chi}) gives us
$\chi_{W(\lambda_i)}$ for $i=1,\dots,l$.
\end{remark}

As we mentioned in the Introduction, such formulas have been also
studied in \cite{DKKMM}, \cite{KKMM1}, \cite{KKMM2}, \cite{KNS},
\cite{T}, \cite{FS1}, \cite{FS2}, \cite{G}, \cite{AKS}, \cite{Cal1}
and \cite{FFJMM}.

\vspace{.3in}

\noindent {\small \sc Department of Mathematics, Ohio State University,
Columbus, OH 43210} \\
{\em Current e--mail address}: corina.calinescu@yale.edu\\
\vspace{.1in}

\noindent {\small \sc Department of Mathematics, Rutgers University,
Piscataway, NJ 08854} \\ {\em E--mail address}:
lepowsky@math.rutgers.edu \\
\vspace{.1in}

\noindent {\small \sc Department of Mathematics and Statistics,
University at Albany (SUNY), Albany, NY 12222} \\ {\em E--mail
address}: amilas@math.albany.edu

\end{document}